\documentclass[11pt]{article}
\usepackage{a4}
\usepackage{graphicx}
\usepackage{parskip}

\usepackage{amsfonts}
\usepackage{natbib}
\usepackage{algorithm}
\usepackage{algorithmic}

\usepackage{latexsym}

\newtheorem{theorem}{Theorem}[section]

\newtheorem{lemma}{Lemma}[section]
\newtheorem{corollary}{Corollary}[section]
\newtheorem{remark}{Remark}[section]

\newtheorem{condition}{Condition}[section]
\newtheorem{definition}{Definition}[section]

\newcommand{\R}{\mathbb{R}}
\newcommand{\Nat}{\mathbb{N}}

\newcommand{\argmin}{\mathop{\arg \min}\limits}

\newcommand{\PP} {{  \rm I\hskip-0.22em P}}

\newcommand{\EE} {{\rm I\hskip-0.48em E}}

\usepackage{a4}

\usepackage{graphicx}

\usepackage{parskip}

\usepackage{amsfonts}

\usepackage{natbib}

\usepackage{algorithm}

\usepackage{algorithmic}

\usepackage{amsmath, amssymb}

\usepackage{latexsym}

\begin{document}
\centerline{\Large On tight bounds for the Lasso}

\vskip .1in
\centerline{Sara van de Geer\footnote{Research supported by
Isaac Newton Institute for Mathematical Sciences,
program {\it Statistical Scalability}, EPSRC Grant Number LNAG/036 RG91310.}}

\centerline{April 3, 2018}

{\bf Abstract} We present upper and lower bounds for
the prediction error of the Lasso.   
For the case of random Gaussian design, we show that under mild
conditions 
the prediction error of the Lasso is up to smaller order terms
dominated by the prediction error of its noiseless counterpart.
We then provide exact expressions for the prediction error
of the latter, in terms of compatibility constants.
Here, we assume the active components of the underlying
regression function satisfy some ``betamin" condition.
For the case of fixed design, we provide upper and lower bounds,
again in terms of compatibility constants.
As an example, we give an up  to a logarithmic term tight bound for the
least squares estimator with total variation penalty. 



\section{Introduction}\label{introduction.section}
Let $X \in \R^{n \times p}$ be an input matrix and $\beta^0 \in \R^p$ a vector of unknown
coefficients. Consider an $n$-vector of noisy observations
$$ Y = X \beta^0 + \epsilon $$
where  the noise $\epsilon \in \R^n $
is a vector of i.i.d.\ standard Gaussians independent of $X$.
The Lasso estimator $\hat \beta$ is 
\begin{equation}\label{Lasso.equation}
\hat \beta \in \argmin_{b \in \R^p}  \biggl \{ \|Y - X b \|_2^2 + 2 \lambda \| b \|_1 \biggr \} 
 \end{equation}
with $\lambda >0$ a regularization parameter (\cite{tibs96}).
Its prediction error is $\| X ( \hat \beta - \beta^0) \|_2^2 $. Main aim of this paper is 
to provide lower bounds for this prediction error, bounds which show
that compatibility constants necessarily enter into the picture.

The results of this paper can be summarized as follows. Firstly, 
suppose the design is random and that
$\Sigma_0 := \EE X^T X/ n$ exists.
Let $\beta^*$ be the noiseless Lasso for random design
\begin{equation}\label{noiselessrandomdesign.equation}
 \beta^* \in \argmin_{b \in \R^p}  \biggl \{ n \|\Sigma_0^{1/2} ( b - \beta^0)\|_2^2 + 2 \lambda \| b \|_1 \biggr \} .
 \end{equation}
For the case where the rows of $X$ are i.i.d\ ${\cal N} ( 0, \Sigma_0)$,
we show in Theorem \ref{concentration2.theorem} that $\| X (\hat \beta - \beta^0) \|_2$
is up to lower order terms equal to $\sqrt n\| \Sigma_0^{1/2} ( \beta^* - \beta^0) \|_2$.
This result is true under the condition that (after normalizing
the co-variance matrix $\Sigma_0$ to having bounded entries)
 the largest eigenvalue
$\Lambda_{\rm max}^2$ of $\Sigma_0$ is of small order $\log n$, and under some
mild condition on the growth of the compatibility constants as $n$ increases.
Secondly, we provide in Theorem
\ref{noiseless.theorem} exact expressions for the prediction error of the noiseless Lasso
in terms of compatibility constants. We require here ``betamin" conditions, which roughly
say that the non-zero coefficients of $\beta^0$ should have the appropriate signs and remain
above the noise level in absolute value. Thirdly, for the case of fixed design,
we present upper and lower bounds for the prediction error $\| X (\hat \beta- \beta^0) \|_2^2$
in terms of weighted compatibility constants. Theorem \ref{noisy.theorem} states the lower bounds,
assuming again certain betamin conditions. The upper bounds
 we present are similar to those obtained the literature and presented
 for completeness. They  are in Corollary \ref{nolambda*.corollary}.
As an illustration we consider least squares estimation
with a (one-dimensional) total variation penalty. We arrive in Corollary \ref{TVupper.corollary}
 at lower and
upper bounds that are the same up to a logarithmic term. 

There are general upper bounds in the literature, in particular {\it sharp oracle bounds}
as in \cite{koltchinskii2011nuclear} (see also 
\cite{giraud2014introduction}, Theorem 4.1 or \cite{vdG2016}, Theorem 2.2).
The oracle bounds involve a compatibility constant, and an 
improved version of this constant
has been developed  in \cite{sunzhang11}, 
\cite{belloni2014pivotal} and \cite{dalalyan2017prediction}.

Main theme of this paper is to gain further insight into the role of the compatibility constant when
applying the Lasso and to see how it occurs in lower bounds. 
In \cite{zhang2014lower} it is shown that for a given sparsity level, there is a design and a lower bound
for the mean prediction error  in the noisy case, that holds for any polynomial time algorithm. This lower bound is close to the known upper bounds and
in particular shows that compatibility conditions or restricted eigenvalue conditions cannot be avoided.
This has also been shown by \cite{bellec2017optimistic}, where a choice of the particular vector of regression
coefficients $\beta^0$ leads to a lower bound matching the upper bound. We further elaborate on this issue,
and provide lower bounds that hold for a large class of vectors $\beta^0$. 

To get an idea of the flavour of the type of bounds we are after, we present in Theorem
\ref{tight.theorem} the case of random design. Details of its proof can be found in Subsection
\ref{tight.section}.
We provide more explicit statements in Theorem \ref{concentration2.theorem}.

Throughout the paper, the active set of $\beta^0$ is denoted by $S_0 := \{ j : \ \beta_j^0 \not= 0 \}$.
Its size is denoted by $s_0 := | S_0| $. Our betamin condition is as follows (its meaning should become
more clear after looking at Section \ref{compatibility.section} where compatibility constants are defined).

\begin{condition} \label{betamin.condition} 
Let
$$ b^* \in \argmin \biggl \{ \| \Sigma_0^{1/2}  b\|_2: \ \sum_{j \in S_0 } | b_j | - \sum_{j \notin S_0} |b_j | =1 \biggr \} $$
   and for $j \in S_0$ let $z_j^* $ be the sign of $b_j^*$. We say that $\beta^0$ satisfies the
   betamin condition for the noiseless case with random design if
   $$ z_j^* \beta_j^0 > { z_j^* b_j^* \over
   \| \Sigma_0^{1/2} b \|_2^2 } { \lambda \over n} , \ \forall  j \in S_0 . $$
\end{condition}

\begin{theorem} \label{tight.theorem} Let the rows of $X$ be i.i.d. \ ${\cal N} (0, \Sigma_0)$, let
$\| \Sigma_0 \|_{\infty}$ be the maximal entry in the co-variance matrix $\Sigma_0$ and
$\Lambda_{\rm max}^2 $ be its largest eigenvalue. For $S \subset \{ 1 , \ldots , p \}$, let
$\kappa^2 (S)$ be the compatibility constant defined in Definition \ref{theoreticalcomp.definition}.
Suppose that
$$ \Lambda_{\rm max}^2 / \| \Sigma_0 \|_{\infty} = o (\log (2p) ) , $$
and
$$ \max \biggl \{ \biggl ( {\| \Sigma_0 \|_{\infty}\over\kappa^2 (S)}  \biggr ) {\log (2p)  |S|   \over
  n  }   
 : \ S \subset \{ 1 , \ldots , p \} , \ |S| \le \biggl ({  \Lambda_{\rm max}^2  \over  \kappa^2 (S_0) } \biggr )  4s_0
  \biggr \}  = o(1).$$
 For some $t>0$, take the tuning parameter $\lambda$ to satisfy
  $$3 \| \Sigma_0 \|_{\infty}^{1/2} \biggl (  \sqrt {2n (\log (2p) + t )} + 2(\log (2p) +t) \biggr )\le \lambda =
  {\mathcal O} \biggl (\sqrt {\| \Sigma_0 \|_{\infty}^{1/2} \log (2p) }\biggr  ) 
   .$$ 
     Then, under Condition \ref{betamin.condition} (the betamin condition for the
     noiseless case with random design), we have
  $$ \| X ( \hat \beta - \beta^0) \|_2^2 = { \lambda^2/n \over \| \Sigma_0^{1/2} b^* \|_2^2  } 
  (1+ o_{\PP} (1)) + {\mathcal O}_{\PP} (1)  $$
  (where in fact $s_0\| \Sigma_0^{1/2} b^* \|_2^2 = \kappa^2 (S_0) $). 
  \end{theorem}
  
  \section{Organization of the paper}
  In Section \ref{compatibility.section} the definition of compatibility constants
  is given and also some of their properties are discussed.
  Section \ref{randomdesign.section} shows that for the case
  of random design the squared ``bias" of the Lasso dominates its
  ``variance". Section \ref{noiseless.section} then gives expressions for
  this ``bias", i.e.\ for the noiseless Lasso. Here, we examine fixed design but the results 
  carry over immediately to random design. In Section \ref{noiselessTV.section}
  the result of Section \ref{noiseless.section} is illustrated with the total
  variation penalty (in one dimension). Section \ref{noisy.section} presents lower
  bounds for the case of fixed design, and Section \ref{compare.section}
  presents some upper bounds. Corollary \ref{nolambda*.corollary}
  is essentially as in the papers
  \cite{sunzhang11}, 
\cite{belloni2014pivotal} and \cite{dalalyan2017prediction},
albeit that do not consider the approximately sparse case.
Section \ref{noisyTV.section} has upper and lower bounds for the least squares
estimator with total variation penalty in the noisy case. Section
\ref{conclusion.section} concludes.
Section \ref{proofs.section} contains the proofs. 

\section{Compatibility constants}\label{compatibility.section}
We introduce some notation in order to be able to define the compatibility constants. This notation
will also be helpful at other places. 
For $S \subset \{1 , \ldots , p \}$  and a vector $b \in \R^p$ let $b_S \in \R^p$ be the vector
with entries
$ b_{j,S} := b_j {\rm l} \{ j \in S \} $, $j=1 , \ldots , p$. We apply the same notation for the $|S|$-dimensional vector 
$\{ b_j \}_{j \in S}$. We moreover write $b_{-S} := b_{S^c}$ where $S^c$ is the the complement
of the set $S$. 

\subsection{Theoretical compatibility constants}

The population version of the compatibility constant will be used
for the case of random design $X$. We call the population
 version the theoretical compatibility constant.

 \begin{definition} \label{theoreticalcomp.definition}
  Let
$\Sigma_0 := \EE X^T X / n$ (assumed to exist).
 Let $S \subset \{ 1 , \ldots , p \}$ be a set of indices and $u \ge 0$ be a constant.
 The theoretical compatibility constant is
 $$\kappa^2 (u,S) := \min \biggl \{  |S| \| \Sigma_0^{1/2} b \|_2^2 : \ \| b_S \|_1 -u \| b_{-S} \|_1 =1 \biggr \} .$$
 For $u=1$ we write $\kappa (1, S) =: \kappa (S)$. 
\end{definition}

\subsection{Empirical compatibility constants} 

For a vector $w $ we let $W:= {\rm diag} (w) $ be the diagonal
matrix with $w$ on the diagonal. 

\begin{definition} \label{compatibility-improved.definition} (\cite{belloni2014pivotal}, \cite{dalalyan2017prediction})
Let $S
\subset \{ 1 , \ldots , p \}$ be a set of indices and $w\in \R^{p- |S|}$  be a vector of non-negative weights. 
The (empirical) compatibility constant is
 is
$$ \hat \kappa^2 (w, S) := \min \biggl \{ |S| \| X b\|_2^2 /n: \  \| b_S \|_1   - 
 \|W b_{-S} \|_1 =1 \biggr \} .$$
 For the case where $w={\bf 1}$ where ${\bf 1}$ denotes
 a vector with all entries equal to one,  put $\hat \kappa^2 (S):= \hat \kappa^2 ({\bf 1}, S)$. 
\end{definition}

\subsection{Some properties of compatibility constants}
One readily sees that the theoretical and empirical compatibility constants differ only in terms
of the matrix used in the quadratic form (which is $\Sigma_0$ in the theoretical
case and the Gram matrix $\hat \Sigma := X^T X / n$ in the empirical case). Thus, when discussing their basic properties
it suffices to deal with only one of the two. In this section,
we therefore restrict attention to the empirical version $\hat \kappa (w,S)$.
Note that we have generalized the empirical version as compared
to the theoretical one, by considering general weight vectors, not just constant
vectors. 
With some abuse of notation, we write $\hat \kappa (u, S)= \hat \kappa ( u {\bf 1} , S) $
when the weights are the constant vector $u {\bf 1} $ (it should be clear from
the context what is meant). 

The empirical compatibility constant as given in Definition \ref{compatibility-improved.definition}
is from \cite{belloni2014pivotal} or \cite{dalalyan2017prediction}.
Another version, from  for instance \cite{vandeGeer:07a} or  \cite{vdG2016} and its references,
is presented in the next definition.

\begin{definition}\label{compatibility.definition} Let $S \subset \{ 1 , \ldots , p \}$ be a set of indices and
$u >0 $ be a constant.
The (older) compatibility constant is
$$ \hat \phi^2 (u,S) := \min \biggl \{ |S| \| X b \|_2^2/n : \ \| b_S \|_1 =1 , \ 
\| b_{-S} \|_1 \le 1/u  \biggr \} .$$
Let $\hat \phi^2 (S) := \hat \phi^2 (1, S)$ be the compatibility constant for the case $u=1$.
\end{definition} 

The constant $\hat \phi (u,S)$ compares, for $b$'s
satisfying a ``cone condition" $\| b_{-S} \|_1 \le  \| b_S \|_1 /u$, the $\ell_2$-norm
$\| X b \|_2$ with the $\ell_1$-norm $\| b_S \|_1$. The constant $\hat \kappa (u, S)$ is similar,
but takes in the comparison more advantage of a ``cone condition"
$\| b_S \|_1 - u \| b_{-S} \|_1 >0$. 
When $\hat \kappa^2 (S)>0$  the null space property holds  
(\cite{donoho2005neighborliness}). We will need throughout that the
compatibility constant is strictly positive at $S_0$
(if it is zero our results cease
to be of any interest). 
This means that we implicitly require throughout 
\begin{condition}\label{inverse.condition} The matrix
   of $X_{S_0}^T X_{S_0}$ is invertible.
    \end{condition}
    
   Here, for any $S \subset
\{ 1 , \ldots , p \}$ the matrix $X_S = \{ X_j \}_{j\in S}$ is the $n \times |S|$ matrix consisting of the columns
of $X$ 
corresponding to the set $S$.

The newer version $\hat \kappa (u, S)$ is an improvement over $\hat \phi(u,S)$
in the sense that $\hat \kappa (u,S)$ is the larger of the two.

\begin{lemma} \label{improve.lemma} For all $u>0$ it is true that
$$ \hat \kappa^2 (u, S) \ge \hat \phi^2 (u, S) . $$
\end{lemma} 

Let now for some $v>0$
$$ b^* \in \arg \min \biggl \{ \| X b \|_2^2 : \ \| b_S \|_1 -  v\| b_{-S} \|_1 =1 \biggr \} . $$
Then by definition
$$ \hat \kappa^2 (v, S) =|S|  \| X b^* \|_2^2 / n .$$
The restriction $\| b_S \|_1 -  v\| b_{-S} \|_1 =1$ does not put any
bound on the $\ell_1$-norm of $b_S^*$. However, if there is a little room to spare,
its $\ell_1$-norm is bounded. This will be useful to understand the betamin conditions
(Conditions \ref{betamin.condition} and \ref{betaminnoiseless.condition}). 
For simplicity we examine only the value $v=1$. 

\begin{lemma} \label{normalize.lemma} Let
$$ b^* \in \argmin \biggl \{ \| X b \|_2^2 /n :\ \| b_S \|_1 - \| b_{-S} \|_1 =1 \biggr \} . $$
Then for $0\le u < 1$
$$ \| b_S^* \|_1 \le { \hat \kappa (S) - u \hat \kappa (u,S) \over(1-u)  \hat \kappa (u,S)} . $$
\end{lemma}

\subsection{Comparing  empirical and theoretical and compatibility}
 Having random quadratic forms in mind, the fact that $\| b_S \|_1 -  \| b_{-S} \|_1 =1$
 gives no bound on the $\ell_1$-norm can be a problem.
Again, if there is a little room to spare in the value of $u $ in the 
compatibility constant,
one {\it does} get a bound on the $\ell_1$-norm. 
We show this in Lemma \ref{kappabound.lemma}, and with this tool in hand we  lower bound
the empirical compatibility constant in terms of the theoretical one
in Lemma \ref{hatkappa.lemma}. 

\begin{lemma}\label{kappabound.lemma} Let $v >u >0$. Then
$$\hat \kappa^2 ( v , S) \ge \min \biggl \{ |S| \| X b \|_2^2 /n : \
\| b_S\|_1- u \| b_{-S} \|_1= 1, \ 
 \| b \|_1 \le 1+ (1+u) / (v-u)   \biggr \}  .$$
\end{lemma}

The following lemma will be applied when bounding
the prediction error of $\hat \beta$  in terms of that of the noiseless Lasso $\beta^*$.
The lemma may also be of interest in itself with applications elsewhere.

\begin{lemma} \label{hatkappa.lemma} Suppose the rows
of $X$ are i.i.d.\ ${\cal N} (0, \Sigma_0) $. Let $\|\Sigma_0 \|_{\infty}$ be the
largest entry in the matrix $\Sigma_0$.
For $v > u$, $(1+u) / (v-u)= {\mathcal O} (1)$
and 
$$ \biggl (    {\| \Sigma_0 \|_{\infty}   \over \kappa^2 ( u, S) } \biggr ) {s\log (2p)\over n}
= o(1), $$
it is true with probability tending to one
that
$$ \hat \kappa^2 ( v  , S) \ge (1- \eta)^2 \kappa^2 ( u {\bf } , S) . $$
where $\eta= o(1)$.
\end{lemma}

 \section{Comparison with the noiseless Lasso when the design is random}
 \label{randomdesign.section}
 In this section we assume that the rows of $X$ are i.i.d.\ copies
 of a Gaussian row vector with mean zero and co-variance matrix
 $\Sigma_0$. We denote the  largest eigenvalue of $\Sigma_0$ by
 $\Lambda_{\rm max}^2 $ and let $\| \Sigma_0 \|_{\infty}$ be its largest entry. We define 
 a noiseless version $\beta^*$ of the Lasso where also the random design is replaced by
 its population counterpart:
 $$ \beta^* \in \argmin_{b \in \R^p} \biggl \{ n \| \Sigma_0^{1/2} ( b- \beta^0)\|_2^2 + 2 \lambda \| b \|_1 \biggr \} . $$ 
 The normalization with $n$ is
 to put things on the scale of the empirical version, 
 as $\EE X^T X = n \Sigma_0$. 
 One may think of $\| X ( \beta^* - \beta^0) \|_2 $ as  ``bias" and
  $\| X ( \hat \beta - \beta^* ) \|_2^2$ as ``variance". We first investigate in some
  detail the ``variance" part in Theorems \ref{randomdesign.theorem} and
  \ref{concentration.theorem}. Then we apply the triangle inequality
  as a way to establish that the squared ``bias" dominates the ``variance", 
  see Theorem \ref{concentration2.theorem}.

   \begin{theorem}\label{randomdesign.theorem} 
   Suppose that
  $$ \rho^2 := \max \biggl \{ \biggl ( {\| \Sigma_0 \|_{\infty}    \over
  \kappa^2 (S) }\biggr )  {  \log (2p)  |S| \over n}
 : \ S \subset \{ 1 , \ldots , p \} , \ |S| \le \biggl ( {  \Lambda_{\rm max}^2  \over  \kappa^2 (S_0) } \biggr ) 4s_0
  \biggr \}  = o(1).$$
 Take for some $t>0$
  $$\lambda \ge  3 \| \Sigma_0 \|_{\infty}^{1/2} \biggl (  \sqrt {2n (\log (2p) + t )} + 2(\log (2p) +t) \biggr ) $$ 
  and define
  $$ \gamma:= (2 \Lambda_{\rm max} ) \sqrt n / \lambda +(2/  \| \Sigma_0 \|_{\infty}^{1/2})
 \rho \lambda / \sqrt {n \log (2p) } . $$
  Then 
 we have for all $x >0$ with probability at least $1- 4 \exp[-t] - \exp[-x] - o(1)$ that
  $$ \| X (\hat \beta - \beta^* ) \|_2 \le \gamma 
    \sqrt {n} \| \Sigma_0^{1/2} (\beta^* - \beta^0) \|_2   + \sqrt {2x} .$$
  \end{theorem}

   Using concentration of measure, one can remove the dependency of the
   confidence level on the value of $t$. This value appears in the choice of the tuning parameter 
   $\lambda$. We make some rather arbitrary choices for the constants.
   
  \begin{theorem}\label{concentration.theorem}
  With the conditions and notations of Theorem
  \ref{randomdesign.theorem}, and assuming in addition that 
  $4\exp[-t] < 1/8 $ (say),
  for $n$ large enough and for all $x >0$, with probability at least $1- 2\exp[-x]$,
  $$ \| X (\hat \beta - \beta^* ) \|_2 \le 
   \gamma 
    \sqrt {n} \| \Sigma_0^{1/2} (\beta^* - \beta^0) \|_2   + 4\sqrt {\log 2} + \sqrt {2x} . $$
   \end{theorem}
  
  We can now make a type of bias-variance decomposition.
  The triangle inequality tells us that
  $$\biggl |  \| X ( \hat \beta - \beta^0) \|_2- \| X ( \beta^* - \beta^0) \|_2 \biggr | 
  \le \| X ( \hat \beta - \beta^* ) \|_2 . $$
  We then approximate the empirical ``bias" $\| X ( \beta^* - \beta^0) \|_2$ by 
  the theoretical ``bias" $\sqrt n \| \Sigma_0^{1/2} (\beta^* - \beta_0)\|_2$
  (which is easy as $\beta^*$ and $\beta^0$ are non-random vectors),
  and use Theorem \ref{randomdesign.theorem} or 
  \ref{concentration.theorem} to bound the ``variance" $\| X ( \hat \beta - \beta^* ) \|_2^2$.
  
  \begin{theorem} \label{concentration2.theorem}
 With the conditions and notations of Theorem \ref{concentration.theorem}, we have for
  $n$ sufficiently large, for all $x>0$ with probability at least $1-2\exp[-x]$
  $$
  \biggl | \| X ( \hat \beta - \beta^0 ) \|_2 - \sqrt n \| \Sigma_0^{1/2} ( \beta^* - \beta^0)  \|_2 \biggr | $$
  $$ \le (\gamma + o(1)) \sqrt n \|\Sigma_0^{1/2} (\beta^* - \beta^0) \|_2  + 4 \sqrt {\log2} + \sqrt {2x}. $$
   \end{theorem}

       \begin{corollary} \label{concentration.corollary}
       Recall that we defined $\gamma$ as
    $$ \gamma:= (2 \Lambda_{\rm max} )\sqrt n / \lambda +(2 /\| \Sigma_0\|_{\infty}^{1/2} )
 \rho \lambda / \sqrt { n \log (2p) } . $$
Therefore, with the conditions and notations of Theorem
  \ref{concentration2.theorem}, and assuming  in addition\\ 
  - $ { \Lambda_{\rm max}^2 / \| \Sigma_0\|_{\infty} } = o(\log (2p ))$,\\
  and\\
 - $ \lambda = o ( \sqrt { \| \Sigma_0 \|_{\infty} n \log (2p) } )/\rho $,\\
 we get with probability at least $1- 2\exp[-x] $
 $$ \biggl |  \| X (\hat \beta - \beta^0 ) \|_2 - \sqrt n \| \Sigma_0^{1/2} (\beta^* - \beta^0) \|_2 \biggr | 
 = o( \sqrt {n} \| \Sigma_0^{1/2} (\beta^* - \beta^0) \|_2 ) + 4 \sqrt {\log 2} + \sqrt {2x} . $$
 In words: the squared ``bias" dominates the ``variance".
  \end{corollary} 
  
  \begin{remark} With the help of Lemma \ref{quadraticform.lemma}, one may also
  prove bounds for $\sqrt n  \| \Sigma_0 ( \hat \beta - \beta^0 ) \|_2$ to complete those for
  of $\| X (\hat \beta - \beta^0 ) \|_2$. We refrain from doing this here to avoid 
  digressions. 
   \end{remark}
   
   \section{The noiseless case with fixed design}\label{noiseless.section}
   In this section we study fixed design $X$ and the noiseless Lasso
   \begin{equation}\label{Lasso*.equation}
   \beta^* \in  \argmin_{b \in \R^p} \biggl \{ \| X ( b - \beta^0) \|_2^2 + 2 \lambda^* \| b \|_1 \biggr \} .
   \end{equation}
In principle the noiseless Lasso considered here 
differs from (\ref{noiselessrandomdesign.equation}), although one can say
that for fixed design $\hat \Sigma = \EE \hat \Sigma =: \Sigma_0$,
with $\hat \Sigma:= X^T X / n $ being the Gram matrix.
In what follows in this section, we do not use any specific properties of $\hat \Sigma$ and
the theory goes through for any positive semi-definite matrix,  $\Sigma$ say.
In the upcoming illustration on functions of bounded variation, the fixed
design setup is the natural one.

Note that we supplied the tuning parameter $\lambda^*$ with a supscript $*$.
This is because in Theorem \ref{hidden-oracle2.theorem} we consider a case with different tuning parameters
for the noisy and the noiseless case, say $\lambda$ and $\lambda^*$.

The Karush-Kuhn-Tucker (KKT) conditions for the noiseless Lasso read
\begin{equation}\label{KKT-noiseless.equation}
X^T X   (\beta^* - \beta^0) + \lambda^*  \zeta^* =0  , \ \zeta^* \in \partial \| \beta^* \|_1 ,
\end{equation}
where $\partial \| b \|_1 $ denotes the
sub-differential of $b \mapsto \| b \|_1 $:
$$ \partial \| b \|_1 = \biggl \{ z \in \R^p: \ z^T b = \| b \|_1 , \ \| z \|_{\infty} \le 1 \biggr \} . $$

Recall that
$$ \hat \kappa^2 (S) = |S | \| X b^* \|_2^2 /n $$
where
\begin{equation}\label{b*.equation}
 b^* \in \arg \min_{b \in \R^p} \biggl \{ \| X b \|_2: \ \| b_{S}  \|_1 - \| b_{-S} \|_1 =1 \biggr\} .
 \end{equation}
Note that $b^*$ given in (\ref{b*.equation}) is not unique, for example we can flip the signs of $b^*$
(i.e., replace $b^*$ by $-b^*$). 

In Theorem \ref{noiseless.theorem} below we give a tight result for the noiseless case
under the condition that the active coefficients in $\beta^0$ are sufficiently large in absolute value:
Condition \ref{betaminnoiseless.condition}. 
Here sufficiently large depends on the magnitude of the entries of a solution $b^*$ of
(\ref{b*.equation}) with $S= S_0$. Therefore, it is of interest to know how large $b^*$ is.
Lemma \ref{normalize.lemma} 
considers its $\ell_1$-norm, and in view of this lemma we conclude that if there is a little room to spare,
the $\ell_1$-norm of $\| b_S^* \|_1$ is bounded, or - in other words - $\{ b_j^*  |S| \}_{j\in S} $ is
bounded ``on average".

For the next condition it is useful to know that we show
in Lemma \ref{not0.lemma} that for $b^*$ given in (\ref{b*.equation}), each coefficient
$b_j^*$ with $j \in S$ is nonzero (provided $\hat \kappa^2 (S) >0$). 

\begin{condition}\label{betaminnoiseless.condition}
Suppose $\hat \kappa^2 (S_0) >0 $. 
Let $b^*$ satisfy (\ref{b*.equation}) with $S= S_0$.
Denote, for $j \in S_0$, the sign of $b_j^*$ as $z_j^*$.
We say that $\beta^0$ satisfies the betamin condition for the noiseless case with fixed design if
$$z_j^* \beta_j^0 > { z_j^* b_j^* s_0 \over  \hat \kappa^2 (S_0)}
{\lambda^* \over n  } \ \forall \ j \in S_0 . $$
\end{condition}

Here is the main theorem for the noiseless case.

\begin{theorem}\label{noiseless.theorem} 
Suppose $\hat \kappa^2 (S_0) >0$. Let $b^*$ satisfy (\ref{b*.equation}) with $S= S_0$. 
If $\beta^0$ satisfies Condition
\ref{betaminnoiseless.condition} (the betamin condition for the noiseless case with fixed design),
 then there exists a solution $ \beta^*$ of the KKT
conditions (\ref{KKT-noiseless.equation}) such that
$$\| X ( \beta^* - \beta^0 ) \|_2^2 = { s_0 \over \hat \kappa^2 (S_0) }
{  \lambda^{*2}  \over n}  . $$
\end{theorem}

\section{The total variation penalty in the noiseless case} \label{noiselessTV.section}
In this section Theorem \ref{noiseless.theorem} is illustrated with the total variation penalty.
For a vector $f \in \R^n$, its total variation is defined as
$${\rm TV}(f) := \sum_{i=2}^n | f_i - f_{i-1} | . $$
Fix a vector $f^0 \in \R^n$ and let
$f^*\in \R^n $ is the least squares approximation of $f^0$ with total variation penalty
\begin{equation}\label{TVnoiseless.equation}:
  f^* \in \argmin_{f \in \R^n} \biggl \{ \| f - f^0 \|_2^2 + 2 \lambda^* {\rm TV} (f) \biggr \} . 
  \end{equation}
Theorem \ref{TV-compatibility.theorem} presents an explicit expression for the compatibility constant $\hat \kappa^2 (S_0)$
where $S_0$ is the set consisting of the locations of the jumps of $f^0$.
Invoking Theorem \ref{noiseless.theorem} one then
arrives at an 
explicit expression for
$\| f^* - f^0 \|_2^2$ provided the jumps of $f^0$ are sufficiently large,
see Corollary \ref{TVnoiseless.corollary}. 

First, we need to rewrite problem (\ref{TVnoiseless.equation}) as a (noiseless) Lasso problem.
Indeed, for $j=1  \ldots , n$,
$$f_j= \sum_{i=1}^n ( f_i - f_{i-1} ) {\rm l} \{ j \ge i \} =: (X b)_j , 
$$
where $X_{j,i} = {\rm l} \{ j \ge i \}$ and $b_i = f_i- f_{i-1}$,
with $f_0 :=0$. 
Hence we can say that $f^0 = X \beta^0$ and $f^* = X \beta^*$ with
$$ \beta^* := \argmin_{b \in \R^n} \biggl \{ \| X( b - \beta^0) \|_2^2 + 2 \lambda^*
\sum_{i=2}^n | b_i | \biggr \} . $$
Note that the first coefficient $b_1$ is not penalized.
It is therefore typically active, and we consider the active set as
the location of the jumps augmented with the index $\{ 1 \}$. We slightly adjust the definition of the compatibility constant
to deal with the a  coefficient without penalty: we set for
 $S \subset \{ 2, \ldots , n \}$ 
 \begin{equation}\label{compatibilityTV.equation}
  \kappa^2 (S):= \min \biggl \{ { |S \cup \{ 1 \} | \| X b \|_2^2  } : \| b_S \|_1 -
  \| b_{-(S \cup \{1 \} ) } \|_1 = 1 \biggr \} . 
  \end{equation}
  
   Let now $S:= \{ d_1+1 , d_1+d_2+1 , \ldots , d_1 +\cdots + d_s+1 \} $ for some 
  $\{ d_j \}_{j=1}^s \subset \{ 2 , \ldots , n \}$ satisfying $\sum_{j=1}^s d_j +2 < n $. 
  The set $S$ represents locations of jumps, $d_1$ is the location of
  the first jump and $\{ d_j\}_{ j=2}^s$
  are the distances between jumps. 
   Let $d_{s+1} := n- \sum_{j=1}^s d_j $ the distance
   between the last jump and the end point.  For simplicity we assume
   that $d_j$ is even for all $j \in \{ 2 , \ldots , s \} $.

  \begin{theorem}\label{TV-compatibility.theorem}
 The compatibility constant $\hat \kappa^2 (S)$ is, up the constant $4$ and the scaling by $1/n$, the harmonic mean of 
of the distances between jumps, including the distance between starting point and
 first jump and last jump and endpoint:
  $$
 \hat  \kappa^2 (S)   
  =  { s+1 \over { n \over d_1}  + \sum_{j=2}^s { 4n\over d_j } + { n \over d_{s+1} }}.$$
  In fact
  $$
  \hat  \kappa^2 (S) = (s+1) \| X b^* \|_2^2 /n
  $$
  where $b_j^* =0 $ for all $j \notin S$ and $  b^*  = \tilde b/ \| \tilde b \|_1 $ with
  \begin{eqnarray*}
  \tilde b_{d_1+1} &= & {n \over d_1 }+ { 2n \over d_2 } ,\\
  \tilde b_{d_2+1} &=& - \biggl ( {2 n\over d_2 } + {2n \over d_{3} }\biggr )  ,  \\
 & \vdots &  \\
  \tilde b_{d_s} &=& (-1)^{s+1} \biggl (  {2n \over d_s} + {n \over d_{s+1} } \biggr )  .
  \end{eqnarray*}
   \end{theorem} 
 
 \begin{corollary} \label{TVnoiseless.corollary}
  Suppose $f^0$ jumps at $S_0:= S= \{ d_1+1 , d_1+d_2+1 , \ldots , d_1 +\cdots + d_{s}+1 \} $, with $s=s_0$.
  Assume $f^0$ alternates between jumps up and jumps down.
  Suppose moreover that
  \begin{eqnarray*}
 | f_{d_1+1}^0 - f_{d_1}^0 | \ &\ge  &\biggl ( {n \over d_1 }+ { 2 n\over d_2 } \biggr )
 { \lambda^*\over n }  ,\\
  | f_{d_2+1}^0 - f_{d_2}^0| \  &\ge & \biggl ( {2n \over d_2 } + {2 n\over d_{3} } \biggr ){ \lambda^*\over n },  \\
 & \vdots &  \\
| f_{d_{s_0} +1}^0 - f_{d_{s_0}}^0 | &\ge &  \biggl ({2 n\over d_{s_0}} + {n \over d_{s_0+1} } \biggr ) { \lambda^*\over n }  .
  \end{eqnarray*}
  Then by Theorem \ref{noiseless.theorem}
 combined with Theorem \ref{TV-compatibility.theorem}
  $$ \| f^* - f^0 \|_2^2 = \biggl ( 
  { n \over d_1}  + \sum_{j=2}^{s_0} { 4n \over d_j } + { n \over d_{s_0+1} }\biggr )
   { \lambda^{*2} \over n} . $$
  At this point it may be helpful to look how this normalizes. Say we choose
  $\lambda^* = \sqrt { n \log n } $. Suppose $\max_{1 \le j \le s_0+1 } {n/ d_j} = {\mathcal O} (s_0+1)  $.
  Then the jumps of $f^0$ are required to be of order at least
  $(s_0+1) \sqrt {\log n/n }  $. We then obtain
  $$ \| f^* - f^0 \|_2^2  = {\mathcal O} \biggl ( (s_0+1)^2 \log n  \biggr ) . $$
  \end{corollary}

      \section{A lower bound in the noisy case with fixed design}\label{noisy.section}
      
   We now turn to the Lasso $\hat \beta $ in the noisy case, given by
      $$ \hat \beta \in \argmin_{b \in \R^p }\biggl  \{ \| Y - X b \|_2^2 + 2 \lambda \| b \|_1 \biggr \} $$
   where 
   $$ Y = X \beta^0 + \epsilon . $$
   We investigate the case of fixed design $X$.
  Recall that we assume throughout  i.i.d.\ standard Gaussian noise.

 \subsection{Towards betamin conditions}\label{betamin.section}
  Consider some vector $\bar v \in \R^{p-s_0}$ with $0< \bar v_j < 1$
  for all $j $. This vector represents the ``noise"  that is to be overruled
  by the penalty.  Define the collection of weights
    $${\cal W} (\bar v ) := \biggl \{ w \in \R^{p- s_0}  : \  \ 1- \bar v_j \le w_j \le 1+ \bar v_j  \  \forall \ j \biggr  \} . $$
     Let for $\bar W := {\rm diag} (1+ \bar v  ) $
   $$ b^* (\bar v)\in \argmin \biggl \{ \| X b \|_2^2 : \ \| b_{S_0} \|_1 - \|  \bar W  b_{-S_0} \|_1 =1 \biggr \} ,
   \  z_j^* (\bar v ):= {\rm sign} (b_j^* (\bar v)  ) , \ j \in S_0 . $$
   Then by definition $ \hat \kappa^2 (1+ \bar v , S_0) = s_0 \| X  b^* (\bar v)  \|_2^2/n$. 
   We remark here that by a slight adjustment of Lemma \ref{not0.lemma}, 
   the assumption $\hat \kappa (1 + \bar v, S_0) >0$ ensures that $b_j^* (\bar v ) \not= 0$ for all $j \in S_0$.
   
   For $w  \in {\cal W} (\bar v) $ we define the convex problem with linear and convex constraints
   $$b(w ) \in \argmin \biggl \{ \| X b \|_2^2: \ z_{S_0}^{*T} (\bar v)  b_{S_0} - \| W b_{-S_0} \|_1 \ge 1 \biggr \} . $$
   Finally, define 
   $$ {\bf b}_j  (\bar v )  := \max_{w \in {\cal W}(\bar v)  } |b_j (w) | / \| X b (w) \|_2^2 , \ j\in S_0 .  $$
   
        \subsection{Projections}\label{projections2.section}
        We denote the projection of
    $X_{-S_0} $ on the space spanned by the columns of $X_{S_0}$
    by $ X_{-S_0}  {\rm P} X_{S_0} $. The projection is always defined but 
    as it is implicitly assumed that  $X_{S_0}^T X_{S_0}$ is invertible (Condition \ref{inverse.condition}),
    we can clarify what we mean by projection by writing
    $$ X_{-S_0}  {\rm P} X_{S_0} := X_{S_0} (X_{S_0}^T X_{S_0} )^{-1} X_{S_0}^T X_{-S_0}  .$$
    The  anti-projection is denoted by
    $$ X_{-S_0} {\rm A} X_{S_0} = X_{-S_0} - X_{-S_0} {\rm P} X_{S_0} . $$
      We define the matrix
    \begin{eqnarray*}
     V_{-S_0, -S_0}  &:= & \biggl ( X_{-S_0} {\rm A} X_{S_0} \biggr )^T \biggl (  X_{-S_0} {\rm A} X_{S_0} \biggr ) \\
     &= &
    X_{-S_0}^T  \biggl (I- X_{S_0}( X_{S_0}^T X_{S_0} )^{-1} X_{S_0}^T \biggr ) X_{-S_0} , 
    \end{eqnarray*}
    and let $\{ v_j^2 \}_{j \notin S_0} $ be the diagonal elements of this matrix.

   \subsection{A lower bound}
   
   The main result for the noisy case is presented in the next theorem. Here, we use the notations and
   definitions of the previous two subsections.
   
   \begin{theorem} \label{noisy.theorem} 
   Take for some $t>0$,
   \begin{equation}\label{lambda.equation}  
     \lambda>  \| v_{-S_0} \|_{\infty} 
   \sqrt {2  (\log (2p) +t) } . 
   \end{equation}
   Define
   $$ \bar v_j := v_j \sqrt {2  (\log (2p) +t) }/ \lambda , \ j \notin S_0  $$
   and
   $$ \bar u_j : =  u_j \sqrt {2  (\log (2p) +t) }/ \lambda, \ j \in S_0. $$
   where $\{ u_j \}_{j \in S_0}$ are the diagonal elements of the matrix
   $(X_{S_0}^T X_{S_0} )^{-1}$.
   Assume that  $\hat \kappa (1+\bar v , S_0) >0$ and that
the following betamin condition holds:
   $$ |\beta_j^0 | > \lambda ( {\bf b}_j (\bar v)+ \bar u_j ) , \ {\rm sign} (\beta_j^0 ) =  z_j^* (\bar v)  \ \forall j \in S_0 .$$
 Then for all $x >0$ with probability at least $1-  \exp[-t] - \exp[-x]$
   there is a solution $\hat \beta$ of the KKT conditions
   such that
   \begin{equation}\label{noisylowerbound.equation}
    \| X ( \hat \beta - \beta^0 ) \|_2 \ge    \sqrt {s_0  \over
    \hat \kappa^2 (1+\bar v , S_0) } \sqrt { \lambda^2 \over n}  -\sqrt {s_0} - \sqrt {2x}  . 
    \end{equation}
   \end{theorem}
   
   Note that for $j \in S_0$, the quantity $u_j$ is the variance of the ordinary least squares estimator 
   of $\beta_j^0$ for the case $S_0$ is known. Thus the betamin condition of 
   Theorem \ref{noisy.theorem} needs that the magnitude of the active coefficients should exceed
   the noise level of the ordinary least squares estimator for known $S_0$. 
   
      \section{Comparison with the noiseless Lasso when the design is fixed}\label{compare.section}
 
  This section studies the case of fixed design and
compares
 the noisy Lasso
 $$ \hat \beta := \arg \min_{b \in \R^p } \biggl \{
 \| Y - Xb \|_2^2 + 2 \lambda \| b \|_1 \biggr \} $$
 with the noiseless Lasso
 $$ \beta^* := \arg \min_{b\in \R^p}  \biggl \{  \| X (b - \beta^0) \|_2^2 + 
 2 \lambda^* \| b \|_1\biggr \}  $$
 where $\lambda^* \le \lambda$.  
 We let $S_* $ be active set of $\beta^*$ and its cardinality
 $s_* := |S_* | $.
 We investigate the error $\| X ( \hat \beta- \beta^*) \|_2$
 in Theorem \ref{hidden-oracle2.theorem}. 
 For $\lambda^*= 0$ we see that $\beta^*= \beta^0$ and then Theorem
  \ref{hidden-oracle2.theorem} gives a bound for $\| X ( \hat \beta - \beta^0)\|_2$.
  This is elaborated upon in Corollary \ref{nolambda*.corollary}. The case $\lambda^* = \lambda$
  is detailed in Corollary \ref{equallambdas.corollary}.
  The error $\| X ( \hat \beta - \beta^*) \|_2^2$ can then seen as ``variance" and
  $\| X ( \beta^* - \beta^0 ) \|_2$ as ``bias".

 \subsection{Projections}\label{projections3.section}
 We now introduce some notations and definitions similar to the ones in
 Subsections \ref{projections2.section}, now for general $S$ instead of just $S=S_0$.
     The projection of  $X_{-S} $ on the space spanned by the columns of $X_{S} $ is denoted by 
     $X_{-S} {\rm P} X_S $.
Recall that such projections are defined, also if $X_S$ does not have full column rank.
 The anti-projection is 
 $$ X_{-S} {\rm A} X_S := X_{-S} - X_{-S} {\rm P} X_S . $$
 Define the matrix
 $$V_{-S, -S}^S := \biggl (X_{-S} {\rm A} X_S \biggr )^T \biggl ( X_{-S} {\rm A} X_S\biggr )  $$
 and let $\{ (v_j^S)^2 \}_{j \notin S} $ be the diagonal elements of this matrix. 
 
 \subsection{Upper bound} 
 
 Recall the KKT conditions for $\beta^*$ as given in (\ref{KKT-noiseless.equation}), involving
 the vector
 $\zeta^*$ in the sub-differential  $\partial  \| \beta^* \|_1$.

   \begin{theorem}\label{hidden-oracle2.theorem} 
   Fix a set $S $ with cardinality $|S| = s$.  Assume that 
  that for some $t>0$
  \begin{equation} \label{lambdaS.equation}
     \lambda>  \| v_{-S}^S \|_{\infty} 
   \sqrt {2  (\log (2p) +t) } 
   \end{equation}
   and write
   \begin{equation}\label{vS.equation}
   \bar v_j^S := v_j^S \sqrt {2  (\log (2p) +t) }/ \lambda , \ j \notin S .
   \end{equation}
   Suppose that
 $$\lambda^* |\zeta _j^*| /  \lambda < 1- \bar v_j^S \  \forall \ j \notin S .$$
 Define
 $$ \bar w_j^S  := { 1- \bar v_j^S  - \lambda^* |\zeta_j^*| / \lambda  \over
 1- \lambda^* / \lambda }  , \ j \notin S. $$
 We have for all $x$ with probability at least $1-  \exp[-t]- \exp[-x] $
 \begin{equation}\label{noisylowerbound2.equation}
    \| X (\hat \beta - \beta^* ) \|_2 \le  \sqrt {s \over \hat \kappa^2( \bar w^S , S)}
    \sqrt {(\lambda- \lambda^* )^2 \over n } 
     + \sqrt { s} +
   \sqrt {2x } . 
   \end{equation}
     \end{theorem}
%
     
       \begin{corollary}\label{nolambda*.corollary} If we take 
   the tuning parameter $\lambda^*$ of the noiseless Lasso equal to zero, 
   Theorem \ref{hidden-oracle2.theorem} gives 
   the following:
with probability at least $1- \exp[-t]-\exp[-x] $
   $$ \| X (\hat \beta - \beta^0 ) \|_2 \le   \sqrt {s_0/\hat \kappa^2(1- \bar v , S_0) }
   \sqrt {\lambda^2 / n} + \sqrt {s_0}+
   \sqrt {2x} . $$
  This result is comparable to results in  \cite{sunzhang11}, 
   \cite{belloni2014pivotal} and \cite{dalalyan2017prediction},
  albeit that we do not deal with the extension to the approximately sparse case. 
One may check that the the combined conclusions of 
 this corollary with that of Theorem \ref{noisy.theorem} also hold with probability
 at least $1-  \exp[-t] - \exp[-x]$. 
 \end{corollary}

\begin{corollary}\label{equallambdas.corollary}
We can also take $\lambda^* = \lambda$ in Theorem \ref{hidden-oracle2.theorem}.
We then formally put $\bar w_j^S = \infty$ for all $j \notin S$ and
 we put $\hat \kappa ( \bar w) = \infty$ as well.
   Let $S $ with $|S| = s$.   Assume that
   \begin{equation}\label{zeta.equation}
 |\zeta _j^*| < 1- \bar v_j^S \  \forall \ j \notin S 
 \end{equation}
 (this implies $S \supset S_*$).
   We have with probability at least $1- \exp[-t] -\exp[-x]$
   $$ \| X (\hat \beta - \beta^* ) \|_2 \le \sqrt { s} + \sqrt {2x}   . $$
        \end{corollary}
        
        Corollary \ref{equallambdas.corollary} is of interest only when $\sqrt s $ is small enough
       This is the case if $\hat \Sigma := X^T X / n$ has a well behaved maximal eigenvalue 
        $\hat \Lambda_{\rm max}^2$. Indeed, one can show in the same way as in Lemma \ref{Lambdamax.lemma} 
        (where $\hat \Sigma$ is replaced by $\Sigma_0$) that
        $$s \le  \biggl ({  \hat \Lambda_{\rm max}^2   \over  ( 1- \| \bar v^S \|_{\infty} )^2 }\biggr ) 
        {n \over \lambda^2}
        \| X (\beta^* - \beta^0 \|_2^2. $$
        Thus if $ \hat \Lambda_{\rm max}^2 / (\| \hat \Sigma \|_{\infty} (1- \| \bar v^S \|_{\infty})^2) = o(\log (2p))$, 
         then
        $ s = o (\| X( \beta^* - \beta^0 ) \|_2^2 ) $. However, for the case of fixed design, one
        might not want to impose such eigenvalue conditions.  Alternatively, one may want
        to resort to irrepresentable conditions. To this end, fix a set $S \supset S_0$. Let
        for $j \notin S$, the projection of the $j^{\rm th}$
        column $X_j$ on $X_S$ be denoted by
        $$ X_j {\rm P} X_S := X_S \gamma_{S,j} . $$
        Then it is not difficult to see that for $j \notin S$
        $| \zeta_j^* | \le \| \gamma_{S,j} \|_1 $. In other words, a sufficient condition for 
        (\ref{zeta.equation}) to hold is the irrepresentable condition
        $$ \| \gamma_{S,j} \|_1 \le 1- \bar v_j^S, \ \forall j \notin S . $$
        We conclude that under irrepresentable conditions the squared ``bias"
        $ \| X ( \beta^* - \beta^0 ) \|_2^2$ dominates the ``variance" $\| X (\hat \beta - \beta^* ) \|_2 $. 
        
           \section{The total variation penalty in the noisy case} \label{noisyTV.section}
     We continue with the total variation penalty of  Section \ref{noiselessTV.section}, but now
     in a noisy setting:
     $$Y= f^0 + \epsilon , $$
     where $f^0 \in \R^n$ is an unknown vector.  The least squares estimator with total variation penalty is
     \begin{equation}\label{TVnoisy.equation}
  \hat  f \in \argmin_{f \in \R^n} \biggl \{ \| Y-f \|_2^2 + 2 \lambda {\rm TV} (f) \biggr \} . 
  \end{equation}
  
  As has become clear from the previous sections, to assess the prediction error in the noisy
  case one needs to evaluate the compatibility constant $\hat \kappa (w, S)$ with weights $w_j \not=
  1$ for $j \notin S$. 
  For the upper bound on the prediction error, we need lower bounds on $\hat \kappa (w, S)$.
 These are derived in \cite{dalalyan2017prediction}, Proposition 2. We re-derive (and slightly improve)
 their result using a different
 proof (the proof in \cite{dalalyan2017prediction} applies a probabilistic argument). 
 
    Suppose as in Section \ref{noiselessTV.section}
    that the locations of the jumps are $S:= \{ d_1+1 , d_1+d_2+1 , \ldots , d_1 +\cdots + d_s+1 \} $ for some 
  $\{ d_j \}_{j=1}^s \subset \{ 2 , \ldots , n \}$ satisfying $\sum_{j=1}^s d_j +2 < n $. 
   Let $d_{s+1} := n- \sum_{j=1}^s d_j $. Assume again for simplicity 
   that $d_j$ is even for all $j \in \{ 2 , \ldots , s \} $.
  
  \begin{lemma}\label{weighted.lemma} Let $w_1, \ldots , w_n$ be non-negative weights. 
  We have 
  $$ { \sqrt {s+1} \over \hat \kappa(w, S)}  \le \| w \|_{\infty} {  \sqrt {s+1}  \over
  \hat \kappa (S) } 
    + \sqrt {n  \sum_{i=2}^n (w_i - w_{i-1} )^2 } , $$
    where as in Theorem \ref{TV-compatibility.theorem}
    $$ {s+1 \over 
 \hat  \kappa^2 (S)  } 
  =  { n \over t_d}  + \sum_{j=2}^s { 4n \over d_j } + { n \over d_{s+1} }. 
  $$
\end{lemma}

 \begin{corollary}\label{TVupper.corollary} 
   Using the notation of Section \ref{compare.section}
   suppose that
   $\lambda$ satisfies (\ref{lambdaS.equation}) with and let
   $\bar v= \bar v^{S_0}$ be given in (\ref{vS.equation}), both with $S:=S_0$.
   Define $\bar v_i=0$ for all $i \in S_0$.
   We then have with $w_i := 1- \bar v_i$, $j \notin S_0 \cup \{ 1 \} $, $w_1=w_2$ and $w_i=1$, $i\in S_0 $ that
   $$ | w_i - w_{i-1} |  \le | v_i - v_{i-1} |/ \| v \|_{\infty} , \  \ i = \{2 , \ldots , n \} .$$
   In \cite{dalalyan2017prediction} it is shown in their Proposition 3 that    
   $$\sum_{i=2}^n ( v_i - v_{i-1})^2 / \| v \|_{\infty}^2 \le (s_0+1) { \log n }  . $$
  Hence one obtains from Lemma \ref{weighted.lemma} with $S=S_0$, combined with Corollary \ref{nolambda*.corollary}, 
  $$ {\sqrt {s_0 +1 }  \over \hat \kappa (1- \bar v, S_0) } \le 
    \sqrt {{s_0+1} \over \hat \kappa (S_0) }
    + \sqrt {(s_0+1) \log n  \over n}   $$
    where as before
    $$ {s_0+1 \over 
 \hat  \kappa^2 (S_0)  } 
  =  { n \over d_1}  + \sum_{j=2}^{s_0} { 4n \over d_j } + { n \over d_{s_0+1} }. 
  $$
    Thus,  with probability at least $1- \exp[-t]- \exp[-x]$
    $$
     \| \hat f - f^0 \|_2 \le  \lambda \biggl (   \sqrt {(s_0+1) \over n \hat \kappa^2 (S_0) }
    + \sqrt {(s_0+1) \log n  \over n } \biggr )    + \sqrt {s_0} + \sqrt {2x} . $$
    Theorem \ref{TV-compatibility.theorem} implies that 
    $$\hat \kappa (1+ \bar v, S_0) \le 
    \hat \kappa (S_0) . $$
    Recall that for the combined conclusion of Theorem \ref{noisy.theorem}
    and Corollary \ref{nolambda*.corollary}
   we do not have to change the confidence level (which is
   $1- \exp[-t]-\exp[-x]$).
We therefore obtain that 
    if the jumps of $f^0$ are sufficiently large in absolute value,
    as given in Theorem \ref{noisy.theorem}, then with probability at least
    $1 -  \exp[-t]-\exp[-x]$
    \begin{eqnarray*}
   \lambda    \sqrt {s_0+1 \over n \hat \kappa^2 (S_0) } - \sqrt {s_0} - \sqrt {2x}  &\le &
     \| \hat f - f^0 \|_2 
     \le   \lambda \sqrt {s_0+1 \over n \hat \kappa^2 (S_0) }
       + \sqrt {s_0} + \sqrt {2x} \\
    &+ & \lambda \sqrt {(s_0+1) \log n  \over n } .
    \end{eqnarray*} 
    \end{corollary}

\section{Conclusion}\label{conclusion.section} 
This paper establishes that in a sense the squared ``bias" of the Lasso dominates the ``variance".
Moreover, lower bounds for the prediction error are given. These lower often match up to
constants or logarithmic factors the upper bounds, or are in fact tight up to smaller
order terms. The bounds show that compatibility constants necessarily enter into the picture.
The lower bounds require ``betamin" conditions, and - for the case
of random design - also certain sparsity conditions.
It is as yet unclear what can be said when betamin conditions to hold. In combination with this, it would also
be of great interest to know what happens when the regression coefficients are not
(approximately) sparse. As far as we know the question to what extend the Lasso will have large
prediction error when sparseness assumptions are violated (i.e. when the Lasso is used
in a scenario not meant for it) is still open. 

{\bf Acknowledgements:} We thank Rico Zenklusen from the Institute of Operations Research, ETH Z\"urich,
and Hamza Fawzi from Department of Applied Mathematics and Theoretical Physics at the University of Cambridge,
for very helpful discussions.

   \section{Proofs} \label{proofs.section} 
   
   \subsection{Proofs of the lemmas in Section \ref{compatibility.section}}
   
   {\bf Proof of Lemma \ref{improve.lemma}.}
   We have to show that $\hat \kappa^2 (u,S) \ge \hat \phi^2 (u,S)$. 
Write
$$ A:= \biggl \{ b :\ \| b_{-S} \|_1 \le  \| b_{S} \|_1 /u, \ \| b_S \|_1 >0  \biggr \} $$
and
$$ B:= \biggl \{ b :\  \| b_{S} \|_1 - u\|  b_{-S} \|_1 >0  \biggr \}  .$$
Then
$$B \subset A . $$
Thus
\begin{eqnarray*}
 \hat \phi^2 (u,S) &=& \min \biggl \{ { |S| \| X b \|_2^2/n \over \| b_S \|_1^2 } : \ 
b \in A \biggr \}\\
& \le & \min \biggl \{ { |S| \| X b \|_2^2/n \over \| b_S \|_1^2 } : \ 
b \in B \biggr \}\\
& = & \hat \kappa^2 (u , S) . 
\end{eqnarray*}
\hfill $\sqcup \mkern -12mu \sqcap$

   {\bf Proof of Lemma \ref{normalize.lemma}.} This lemma bounds the $\ell_1$-norm
   of the minimizer $b^*$ if there is a little room to spare.
   We have
\begin{eqnarray*}
\| b_S^* \|_1 - u \| b_{-S}^* \|_1 &\le&  \sqrt {| S| /n} \| X b^* \|_2 / \hat \kappa (u,S)  \\
&=&  \hat \kappa (S) / \hat \kappa (u,S) .
\end{eqnarray*}
On the other hand
\begin{eqnarray*}
\| b_S^* \|_1 - u \| b_{-S}^* \|_1&= &\| b_S^* \|_1 -  \| b_{-S}^* \|_1+ (1-u) \| b_{-S}^* \|_1 \\
& = & 1+ (1-u) \| b_{-S}^* \|_1 .
\end{eqnarray*}
Thus
$$ \| b_{-S}^* \|_1 \le {  \hat \kappa (S)- \hat \kappa (u,S)  \over (1-u) \hat \kappa (u,S)}   , $$
yielding
$$ \| b_S^* \|_1 =1+\| b_{-S}^* \|_1 \le { \hat \kappa (S) - u \hat \kappa (u,S) \over(1-u)  \hat \kappa (u,S)}.$$
\hfill $\sqcup \mkern -12mu \sqcap$

{\bf Proof of Lemma \ref{kappabound.lemma}.}
This lemma shows that one has a bound for the $\ell_1$-norm
in the ``cone condition" if there is a little room to spare.
Consider a vector $b \in \R^p$ satisfying
$$ \| b_S\|_1- v \| b_{-S} \|_1 =1 .$$
Since
$$\| b_S\|_1- v \| b_{-S} \|_1 =\| b_S\|_1- u \| b_{-S} \|_1 - (v-u) \| b_{-S} \|_1  $$
we obtain
$$(v-u) \| b_{-S} \|_1 = \| b_S\|_1- u \| b_{-S} \|_1-1 \le  \| b_S\|_1- u \| b_{-S} \|_1. $$
Moreover, clearly
$$ \| b_S\|_1- u \| b_{-S} \|_1= (v-u) \| b_{-S} \|_1 +1\ge 1 . $$
It follows that
$$ \min\biggl  \{ \| Xb \|_2: \  \| b_S\|_1- v \| b_{-S} \|_1 =1 \biggr \} $$ $$ \ge
\min \biggl \{ \| Xb \|_2:\ (v-u) \| b_{-S} \|_1  \le  \| b_S\|_1- u \| b_{-S} \|_1, \ 
\| b_S\|_1- u \| b_{-S} \|_1\ge 1 \biggr \} .$$
Suppose now that for some $c>1$
$$(v-u) \| b_{-S} \|_1  \le  \| b_S\|_1- u \| b_{-S} \|_1, \ 
 \| b_S\|_1- u \| b_{-S} \|_1=c  . $$
Define
$$ \tilde b := b/ c . $$
Then
$$ (v-u) \| \tilde b_{-S} \|_1  \le 1, \ 
 \| \tilde b_S\|_1- u \| \tilde b_{-S} \|_1 =1 . $$
 Moreover
 $$ \| X b \|_2 = c \| X \tilde b \|_2 > \| X \tilde b \|_2. $$
 Therefore
 \begin{eqnarray*}
 & &  \min \biggl \{ \| Xb \|_2:\ (v-u) \| b_{-S} \|_1  \le  \| b_S\|_1- u \| b_{-S} \|_1, \ 
\| b_S\|_1- u \| b_{-S} \|_1\ge 1 \biggr \} \\
&=&
\min \biggr \{ \| Xb \|_2:\ (v-u) \| b_{-S} \|_1  \le 1, \ 
\| b_S\|_1- u \| b_{-S} \|_1= 1 \biggr \} . 
\end{eqnarray*}
But if $(v-u) \| b_{-S} \|_1  \le 1$ and $\| b_S\|_1- u \| b_{-S} \|_1= 1$ we see that
\begin{eqnarray*}
 \| b \|_1& \le & \| b_S \|_1 + \| b_{-S} \|_1 =1 + (1+ u) \| b_{-S} \|_1 \\
 &\le  &
1+ (1+u) / (v-u).
\end{eqnarray*}
\hfill $\sqcup \mkern -12mu \sqcap$

{\bf Proof of Lemma \ref{hatkappa.lemma}.} 
This lemma lower bounds the empirical compatibility constant by
the theoretical one. Here is a proof.
If $\| b_S\|_1- u \| b_{-S} \|_1= 1$ we know that
$$ 1 \le \| \Sigma_0^{1/2} b \|_2 \sqrt {s} / \kappa(u , S ) .$$ 
It therefore follows from Lemma \ref{kappabound.lemma} that
$$\hat \kappa^2 ( v , S) \ge \biggl \{ |S| \| X b \|_2^2 /n : \
\| b_S\|_1- u \| b_{-S} \|_1= 1, \ 
 \| b \|_1 \le M(u,v)
\| \Sigma_0^{1/2} b \|_2  \biggr \}  $$
where
$$M(u,v):=
(1+ (1+u) / (v-u)) \sqrt s / \kappa (u  , S)=
o( \sqrt { n /( \|\Sigma_0 \|_{\infty} \log (2p) )} ) . $$
 In view of Lemma \ref{quadraticform.lemma} we know that when
 $M = o (\sqrt { n / ( \| \Sigma_0 \|_{\infty} \log (2p) } )$, then
with probability tending to one
$$ \inf_{\| b \|_1 \le M\| \Sigma_0^{1/2}  b \|_2} { \| X b \|_2^2/ n   \over 
\| \Sigma_0^{1/2} b \|_2^2 } \ge (1- \eta_M )^2 $$
for suitable  $\eta_M = o(1)$.
Hence with probability tending to one
$$ \min \biggl \{  \| X b \|_2^2 /n : \
\| b_S\|_1- u \| b_{-S} \|_1= 1, \ 
 \| b \|_1 \le M(u,v)
\| \Sigma_0^{1/2} b \|_2  \biggr \}  $$
$$ \ge (1- \eta_{M(u,v)} )^2 \min \biggl \{  \| \Sigma_0^{1/2} b \|_2^2  : \
\| b_S\|_1- u \| b_{-S} \|_1= 1 \biggr \}  = (1- \eta_{M(u,v)} )^2 \kappa^2(u  , S) . $$
\hfill $\sqcup \mkern -12mu \sqcap$

      \subsection{Proof of Theorem \ref{randomdesign.theorem}.} 
     
 The proof is organized as follows. We first present a bound for
 $\| \Sigma_0 ( \beta^*- \beta_0)\|_2$ in Lemma \ref{population.lemma}. This will be used
 to bound later the number of active variables $s_*$ of $\beta^*$,
 or rather some extended version of it involving sub-differential
 calculus, see Lemma \ref{Lambdamax.lemma}.
 We then establish in Lemma \ref{onT.lemma} a deterministic bound assuming
 we are on some subset of the underlying probability space. Then in Lemma \ref{T123.lemma} we show that this subset
 has large probability.
 
 The noiseless Lasso $\beta^*$ given in (\ref{noiselessrandomdesign.equation})
satisfies the KKT conditions
 \begin{equation}\label{populationKKT.equation}
 n  \Sigma_0 ( \beta^* - \beta^0) +\lambda \zeta^*=0 , \ \zeta^* \in \partial \| \beta^* \|_1 ,
  \end{equation}
   where $\partial \| b \|_1 $ is the sub-differential of $b \mapsto  \| b \|_1$:
  $$ \partial \| b \|_1 := \biggl \{ z: \ \| z \|_{\infty}\le 1 , \ z^T b = \| b \|_1 \biggr \} . $$
  This will be used in Lemma \ref{Lambdamax.lemma} and again in Lemma \ref{onT.lemma}.
  In the latter we also invoke the KKT conditions for $\hat \beta$
  \begin{equation}\label{hatKKT.equation}
  X^T X  ( \hat \beta - \beta^0) + \lambda \hat \zeta =X^T \epsilon , \ \hat \zeta \in 
  \partial \| \hat \beta \|_1 .
  \end{equation}

 \subsubsection{ A bound for the number of active variables of $\beta^*$}
 
 First we bound the prediction error of $\beta^*$.
 
 \begin{lemma} \label{population.lemma} Suppose $\kappa^2 (S_0) >0$.
 Then
 $$ n \| \Sigma_0^{1/2} ( \beta^* - \beta^0) \|_2^2 
 \le {  s_0 \over  \kappa^2 (S_0) } { \lambda^2  \over n}. $$
  \end{lemma}
  
  {\bf Proof of Theorem \ref{population.lemma}.} This follows from a slight
  adjustment of Theorem \ref{hidden-oracle2.theorem} in this paper. This is a big detour however, so let us present
  a self-contained proof as well. 
   By the KKT conditions (\ref{populationKKT.equation})
  $$ - ( \beta^* - \beta^0)^T \zeta^* \le \| \beta^0 \|_1 - \| \beta^* \|_1 \le
  \| \beta_{S_0}^* - \beta^0 \|_1 - \| \beta_{-S_0}^* \|_1 . $$
  So if $  \| \Sigma_0^{1/2} ( \beta^* - \beta^0) \|_2^2 >0 $ we obtain by the definition
  of the compatibility constant $\kappa^2(S_0)$ that
  $$n  \| \Sigma_0^{1/2} ( \beta^* - \beta^0) \|_2^2 \le \lambda \sqrt {s_0 }\| \Sigma_0^{1/2} (\beta^* - \beta^0) \|_2
   / \kappa (S_0) . $$
   This yields the result of the lemma. 
      \hfill $\sqcup \mkern -12mu \sqcap$ 
    
    Consider the set  $S_* := \{ \beta_j^* \not=0 \} $ of
 active coefficients of $\beta^*$. We bound the size of this set. In fact
  we look at bound for the size of a potentially larger set, namely 
  the set $ S_* (\nu)  := \{ j: \ | \zeta_j^*| \ge 1- \nu \} $
  where $0 \le \nu < 1$ is arbitrary.  Note that indeed $S_* \subset S_* (\nu) $.
  We pin down the value of $\nu$ to $\nu=1/2$ but the argument goes through for
  other values if one adjusts the constants accordingly.
  We still keep the symbol $\nu$ at places to facilitate tracking the constants. 
  
  \begin{lemma} \label{Lambdamax.lemma} We have that
  $$ | S_* (\nu ) | \le  {  \Lambda_{\rm max}^2  \over (1- \nu)^2  } 
  { n^2 \over \lambda^2} \| \Sigma_0^{1/2} (\beta^*- \beta^0) \|_2^2  \le {\Lambda_{\rm max}^2  \over (1- \nu)^2}
  {s_0 \over \kappa^2 (S_0) } 
  .$$
   \end{lemma}
   
   {\bf Proof of Lemma \ref{Lambdamax.lemma}.}
   Since
   $$ \| \zeta^* \|_2^2 \ge \| \zeta_{S_* (\nu) }^* \|_2^2 \ge (1- \nu)^2 | S_* (\nu ) | $$
   it follows from the KKT conditions (\ref{populationKKT.equation}) that
   $$ (1- \nu)^2 | S_* (\nu) | \le   \| \Sigma_0 ( \beta^* - \beta^0 ) \|_2^2 { n^2 \over  \lambda^2} \le 
   \Lambda_{\rm max}^2 \| \Sigma_0^{1/2}  (\beta^* - \beta^0)\|_2^2 { n^2  \over  \lambda^2 } .$$
  The proof is completed by applying the upper bound of Lemma \ref{population.lemma}
  $$\| \Sigma_0^{1/2}  (\beta^* - \beta^0)\|_2^2 \le {  s_0 \over  \kappa^2 (S_0) } 
  {\lambda^2 \over n^2 } . $$
   \hfill $\sqcup \mkern -12mu \sqcap$
   
    \subsubsection{Projections}
  Let
  $S:= S_* (\nu) $, $s:= |S|$ (where $\nu=1/2$). 
  Set
  $${\bf U} (S) := \| \epsilon {\rm P} X_S \|_2 $$
  where $\epsilon {\rm P} X_S $ is the projection of $\epsilon $ on the space spanned by the
  columns of $X_{S}$.
  Denote the anti-projection of $X_{-S}$ on this space by
  $$X_{-S} {\rm A} X_S := X_{-S} - X_{-S} {\rm P} X_S .$$
  
  \subsubsection{Choice of $\lambda$}
  Recall we take for some $t>0$
  $$ \lambda \ge 3 \| \Sigma_0\|_{\infty}^{1/2} 
 \biggl ( \sqrt {2    ( \log (2p)+t )} +
    2(\log (2p) +t) \biggr ).$$

   \subsubsection{The sets ${\cal T}_1$, ${\cal T}_2$ and ${\cal T}_3$}
   Write 
   $$v_0:=  \| \Sigma_0 \|_{\infty} ^{1/2}\biggl (  \sqrt {2n  (\log (2p ) + t )}+
   2(\log (2p) +t)  \biggr )/\lambda . $$
   We now define a suitable subset of the underlying probability space, on which we
   can derive the searched for inequality.
 This subset will be the intersection of the following sets:
   \begin{eqnarray*}
  {\cal T}_1 &:= & \biggl \{ \| (X_{-S} {\rm A} X_S)^T  \epsilon \|_{\infty} \le \lambda v_0 , \ {\bf U} (S) \le \sqrt {s} + \sqrt {2x} \biggr \}, \\
 {\cal T}_2 & := & \biggl \{ \ \| (X^T X  - n\Sigma_0 ) (\beta^* - \beta^0) \|_{\infty} \le
  \lambda \delta  \biggr \} ,\\
  {\cal T}_3 &:=& \biggl \{ \hat \kappa^2 ( (v- v_0-\delta)/ \delta, S) \ge
  (1- \eta)^2\kappa^2 (S)   \biggr \} ,
 \end{eqnarray*}
where $x >0$ is arbitrary, 
  $\delta := \| \Sigma_0^{1/2}  (\beta^* - \beta^0)\|_2 $,  and where
  $\eta \in (0,1)$ is arbitrary. We pin down $\eta$ to $\eta=1/2$ like we did with $\nu$.
We require that $\nu- v_0  -2 \delta >0$. Since $\nu=1/2$ and
$v_0 \le 1/3$ this is the case for $\delta \le 1/(12)$. In view of 
Lemma \ref{population.lemma}, Theorem \ref{randomdesign.theorem} is about the 
case $\delta=o(1)$, so $\delta \le 1/ (12)$ will be true for $n$ sufficiently large.

  \subsubsection{Deterministic part} 

 \begin{lemma}\label{onT.lemma} On ${\cal T}_1 \cap {\cal T}_2 \cap {\cal T}_3$
 it holds that
 $$\| X (\hat \beta - \beta^*) \|_2 \le\biggl (  {  \Lambda_{\rm max}  \over (1- \nu)  } { \sqrt n \over \lambda} 
 +   \sqrt {s \over  \kappa^2  (S)} {\lambda \over (1- \eta) n} \biggr ) \sqrt n \delta + \sqrt {2x} .$$ 
  \end{lemma}
 
 {\bf Proof of Lemma \ref{onT.lemma}.} 
The KKT conditions  (\ref{populationKKT.equation})
and
(\ref{hatKKT.equation}), for $\beta^*$ and $\hat \beta$ respectively,  are
 $$ X^T X  ( \beta^* - \beta^0 ) + \lambda \zeta^* = Z , $$
 with $Z:= (X^T X  - n \Sigma_0 )(\beta^* - \beta^0) $, and
 $$ X^T X  (\hat \beta - \beta^0) + \lambda \hat \zeta = X^T \epsilon . $$
 So subtracting the first from the second
 $$
X^T X  ( \hat \beta - \beta^*) +\lambda \hat \zeta-\lambda \zeta^*= X^T \epsilon -Z . 
 $$
 Multiplying with $\hat \beta - \beta^* $ yields
 \begin{equation}\label{KKTdifference.equation}
\| X (\hat \beta - \beta^*) \|_2^2  +\lambda 
 (\hat \beta - \beta^*)^T (\hat \zeta-\zeta^*) = (\hat \beta - \beta^*)^T ( X^T \epsilon -
 Z ). 
 \end{equation}
 We write (as in the proof of Theorem \ref{hidden-oracle2.theorem}) with $S:= S_* (\nu) $, $s:= |S|$, 
 $$ X_S \hat b_S := X_S  (\hat \beta_S - \beta_S^*) +
 ( X_{-S} {\rm P} X_S ) \hat \beta_{-S} . $$
 Since $|\zeta_j^* | \le 1- \nu < 1$ for all $j \notin S$, it must be true that
 $\beta_{-S}^*=0$. Therefore
 $$ X ( \hat \beta - \beta^*) = X_S \hat b_S + ( X_{-S} {\rm A} X_S ) \hat \beta_{-S} . $$
 So
 $$ ( \hat \beta - \beta^* )^T X^T \epsilon =
 \hat b_S^T X_S^T \epsilon + \hat \beta_{-S}^T ( X_{-S} {\rm A} X_S )^T \epsilon . $$
 We use that (on ${\cal T}_1$)
 \begin{eqnarray*}
  \hat b_S^T X_S^T \epsilon &\le&  {\bf U} (S)  \| X_S \hat b_S \|_2 \\
  &\le& {\bf U} (S) \| X ( \hat\beta - \beta^* ) \|_2 \\
  &\le & ( \sqrt s + \sqrt {2x} )
 \| X( \hat \beta - \beta^*) \|_2
 \end{eqnarray*}
 and
 $$ \hat \beta_{-S}^T ( X_{-S} {\rm A} X_S )^T \epsilon \le \| \hat \beta_{-S} \|_1 
 \| ( X_{-S} {\rm A} X_S )^T \epsilon \|_{\infty} \le \lambda v_0 \| \hat \beta_{-S} \|_1.$$
 Moreover (on ${\cal T}_2$)
 $$- (\hat \beta- \beta^*)^T Z \le \| \hat \beta - \beta^* \|_1 \| Z \|_{\infty} \le
 \lambda \delta \| \hat \beta- \beta^* \|_1 .  $$
 Then
 \begin{eqnarray*}
 (\hat \beta - \beta^*)^T ( \zeta^* -\hat \zeta ) &= &
 \hat \beta^T \zeta^* - \beta^{*T} \zeta^* + \beta^{*T} \hat \zeta - \hat \beta^T \hat \zeta \\
 &= & \hat \beta^T \zeta^* - \| \beta^*\|_1  + \beta^{*T} \hat \zeta - \| \hat \beta \|_1 \\
 & \le & \| \hat \beta_S \|_1 - \| \beta_S^* \|_1 + \| \beta_S^* \|_1 - \| \hat \beta_S \|_1 \\
 & + & \hat \beta_{-S}^T \zeta_{-S}^*  - \| \hat \beta_S \|_1 \\
 & =  & \hat \beta_{-S}^T \zeta_{-S}^*  - \| \hat \beta_S \|_1 \\
 & \le & (1- \nu) \| \hat \beta_{-S} \|_1 - \| \hat \beta_{-S} \|_1 \\
 & = & - \nu \| \hat \beta_{-S} \|_1 . 
 \end{eqnarray*}
 
 Inserting these bounds in (\ref{KKTdifference.equation}) gives
 $$\| X ( \hat \beta - \beta^* ) \|_2^2 +  \lambda (\nu - v_0 - \delta)  \| \hat \beta_{-S} \|_1 \le
 (\sqrt s+ \sqrt {2x} ) \| X ( \hat \beta - \beta^* ) \|_2 + \lambda \delta \| \hat \beta_S - \beta_S^* \|_1 . $$
 If
 $$\| X ( \hat \beta - \beta^* ) \|_2  \le  (\sqrt s+ \sqrt {2x} ) $$
 we are done as 
 by Lemma \ref{Lambdamax.lemma},
  $ \sqrt {s} \le  \Lambda_{\rm max} \delta n  / ((1- \nu)   
   \lambda)  $. 
   If 
 $$\| X ( \hat \beta - \beta^* ) \|_2  >  (\sqrt s+ \sqrt {2x} ) $$
 we get
 $$ (\nu -v_0 - \delta )  \| \hat \beta_{-S} \|_1 < \delta \| \hat \beta_S - \beta_S^* \|_1  $$
 or
 $$\| \hat \beta_{S} - \beta^*  \|_1 - (( \nu - v_0 - \delta ) / \delta) \|\hat \beta_{-S} \|_1 >0 . $$
 But (on
 ${\cal T}_3$)
   \begin{eqnarray*}
 & & \| \hat \beta_S - \beta_S^* \|_1 - (( \nu - v_0 - \delta ) / \delta) \| \hat \beta_{-S} \|_1\\
  & \le & 
  { \sqrt {s} \| X (\hat \beta - \beta^*) \|_2 \over \sqrt n \hat \kappa(( \nu - v_0 - \delta ) / \delta) , S)} \\ & \le &
 { \sqrt {s}  \| X (\hat \beta - \beta^*) \|_2 \over  \sqrt n \kappa(S)(1- \eta)}  .
 \end{eqnarray*} 

 This gives
 $$\| X (\hat \beta - \beta^*) \|_2 \le \sqrt {s} + \sqrt {2x} +
 \lambda \delta \sqrt {s} / (\sqrt n \kappa (S)(1-\eta))  .$$
 Again, by Lemma \ref{Lambdamax.lemma},
  $ \sqrt {s} \le   \Lambda_{\rm max} \delta n  / ((1- \nu)   
   \lambda)  $.
  We see that
  $$\| X (\hat \beta - \beta^*) \|_2 \le\biggl (  {  \Lambda_{\rm max}  \over (1- \nu)  } { \sqrt n \over \lambda} 
 +  {  \sqrt {s} \over   \kappa (S)(1- \eta)} {\lambda \over (1- \eta) n}  \biggr ) \sqrt n \delta  + \sqrt {2x} .$$
 \hfill $\sqcup \mkern -12mu \sqcap$
 
 \subsubsection{Random part}
 
 We apply the tools of Section \ref{tools.section}.

  \begin{lemma}\label{T123.lemma} It holds that
$$ \PP \biggl ({\cal T}_1 \cap {\cal T}_2 \cap {\cal T}_3 \biggr ) \ge 1- 4 \exp[-t] - \exp[-x] - o(1) . $$
\end{lemma}

{\bf Proof of Lemma \ref{T123.lemma} .}
We first show that $\PP ({\cal T}_1) \ge 1- 2 \exp[-t] - \exp[-x]$. 
One component of this is to show that with probability at least $1- 2 \exp[-t]$
$$ \| (X_{-S} {\rm A} X_S)^T  \epsilon \|_{\infty}   
\le  \lambda v_0. $$
For a square matrix $B$, let ${\rm diag} (B) $ be its diagonal.
By Lemma \ref{Gauss.lemma} 
we know that with probability at least $1- \exp[-t]$
$$  \| (X_{-S} {\rm A} X_S)^T  \epsilon \|_{\infty} \le \| {\rm diag}((X_{-S} {\rm A} X_S)^T (X_{-S} {\rm A} X_S)) \|_{\infty}^{1/2}
\sqrt { 2 ( \log (2p) +t ) } . $$
But
$$\| {\rm diag} ((X_{-S} {\rm A} X_S)^T (X_{-S} {\rm A} X_S) ) \|_{\infty} \le \| {\rm diag} (X^T X ) \|_{\infty} . $$
Moreover in view of Lemma \ref{chi2.lemma}, and using the union bound,
with probability at least $1-\exp[-t] $
$$\biggl |  {\rm diag} (\| X^T X) \|_{\infty}^{1/2} -\sqrt n  \| {\rm diag} (\Sigma_0 ) \|_{\infty}^{1/2} \biggr | \le 
\| \Sigma_0 \|_{\infty}^{1/2}  \sqrt {2( \log (2p) +t ) } . $$
So with probability at least $1- 2\exp [-t]$
$$ \| (X_{-S} {\rm A} X_S)^T  \epsilon \|_{\infty} \le  \| \Sigma_0 \|_{\infty}^{1/2} 
\biggl ( \sqrt {2n (\log (2p) +t ) } +
2( \log (2p) +t )  \biggr )\le  \lambda v_0. $$
The second component is to show that
$$  \PP (  {\bf U} (S) \le \sqrt {s} + \sqrt {2x} ) \le \exp[-x] , $$
but this follows immediately from Lemma \ref{chi2.lemma}. 

Next we show that $\PP ( {\cal T}_2 ) \le 2 \exp[-t] $. Set
$Z:= ( X^TX - n\Sigma_0 ) (\beta^* - \beta^0) $.
Clearly $X (\beta^* - \beta^0)$ is a Gaussian vector with i.i.d. entries with mean zero and variance
  $\| \Sigma_0^{1/2} (\beta^* -\beta^0)\|_2^2 $. 
  Hence, applying Lemma \ref{2dGauss.lemma} with $\sigma_u^2 \le \| \Sigma_0 \|_{\infty} $, 
$\sigma_v^2= \| \Sigma_0^{1/2} (\beta^* - \beta^0) \|_2^2$ 
and using the union bound, we obtain that with probability at least $1- 2 \exp[-t] $
$$ \| Z\|_{\infty} \le 3 \| \Sigma_0 \|_{\infty}^{1/2} \| \Sigma_0^{1/2} (\beta^* - \beta^0) \|_2
 \biggl ( \sqrt {2n  ( \log (2p)+t } +
    \log (2p) +t \biggr ). $$
    
    Finally, the result $\PP ({\cal T}_3 )=1-o(1)$ follows from Lemma \ref{hatkappa.lemma}.
    \hfill $\sqcup \mkern -12mu \sqcap$
    
    \subsubsection{Collecting the pieces}
   Combining  Lemma \ref{onT.lemma} with Lemma \ref{T123.lemma} completes
   the proof of Theorem \ref{randomdesign.theorem}.

 \subsection{Proof of Theorems \ref{concentration.theorem} and \ref{concentration2.theorem}}
  
  We use the concentration of measure, Lemma \ref{concentration.lemma}.

   {\bf Proof of Theorem \ref{concentration.theorem}.}
 Let $m^* := \EE( \| X( \hat \beta - \beta^* ) \|_2 \vert X)$. 
  Then we have  (by Lemma \ref{concentration.lemma}) that with 
  probability at least $1- 1/8 - 3/4 - o(1)$
    $$ \|X (\hat \beta - \beta^*) \|_2 \ge m^* - 2 \sqrt {\log 2} $$
  as well as  (by Theorem \ref{randomdesign.theorem}), 
    $$\|X (\hat \beta - \beta^*) \|_2 \le 
 \gamma
    \sqrt {n} \| \Sigma_0^{1/2} (\beta^* - \beta^0) \|_2   + 2\sqrt {\log 2} . $$
    Thus
    $$m^* \le \gamma
    \sqrt {n} \| \Sigma_0^{1/2} (\beta^* - \beta^0) \|_2   + 4\sqrt {\log 2} .$$
    Applying again Lemma \ref{concentration.lemma} we see that
    $$\PP\biggl  ( \| X (\hat \beta - \beta^* ) \| \ge 
   \gamma
    \sqrt {n} \| \Sigma_0^{1/2} (\beta^* - \beta^0) \|_2   + 4\sqrt {\log 2} + \sqrt {2x} \biggr ) $$
    $$ \le \PP\biggl  ( \| X (\hat \beta - \beta^* ) \| \ge m^* + \sqrt {2x} \biggr ) \le 2\exp[-x] . $$
    \hfill $\sqcup \mkern -12mu \sqcap$
    
     {\bf Proof of Theorem \ref{concentration2.theorem}.} 
   By the triangle inequality
   $$\biggl | \| X (\hat \beta  - \beta^0) \|_2-  \| X ( \beta^* - \beta^0) \|_2\biggr |  \le \| X ( \hat \beta - \beta^*) \|_2 
 .$$
  By Lemma \ref{chi2.lemma}, with
  with probability at least $1- 2/n $
  $$ \biggl | \| X ( \beta^* - \beta^0 ) \|_2 - \sqrt n \| \Sigma_0^{1/2} (\beta^* - \beta^0 ) \|_2 \biggr | 
  \le ( \sqrt {2 \log n} ) \| \Sigma_0^{1/2} ( \beta^* - \beta^0) \|_2  .$$
  So
  with probability at least $1- 4 \exp[-t] -  \exp[-x] - o(1) - 2/n $
  (subtracting the term $2/n$ to follow the argument, as of course it can be
  included in the $o(1)$ term)
  $$\biggl |   \| X (\hat \beta  - \beta^0) \|_2 - \| \Sigma_0^{1/2} ( \beta^* - \beta^0) \|_2 \biggr | 
  \le    (\gamma + \sqrt {2 \log n/n } ) \sqrt n \| \Sigma_0^{1/2} ( \beta^* - \beta^0) \|_2+ \sqrt {2x}  .$$
 Let $m^0 := \EE (\| X ( \hat \beta - \beta^0) \|_2 \vert X).$
  Using the same arguments as in Theorem \ref{concentration.theorem}, we arrive at
  $$ m^0 -2 \sqrt {\log 2}  \le  (1 + \gamma + \sqrt {2 \log n/n } ) \sqrt n \| \Sigma_0^{1/2} ( \beta^* - \beta^0) \|_2+ 2\sqrt {\log 2} $$
  and
 $$ (1 - \gamma - \sqrt {2 \log n/n } ) \sqrt n \| \Sigma_0^{1/2} ( \beta^* - \beta^0) \|_2- 2\sqrt {\log 2} 
  \le m^0 + 2 \sqrt {\log 2} ,$$
  or
  $$ \biggl | m^0 - \sqrt n \| \Sigma_0^{1/2} ( \beta^* - \beta^0) \|_2 \biggr | \le
   (\gamma + \sqrt {2 \log n/n } ) \sqrt n \| \Sigma_0^{1/2} ( \beta^* - \beta^0) \|_2 + 4 \sqrt{\log 2} . $$
   Thus, inserting the triangle inequality
   \begin{eqnarray*}
   & \ & \biggl | \| X (\hat \beta - \beta^0 ) \|_2 - \sqrt n \| \Sigma_0^{1/2} ( \beta^* - \beta^0) \|_2 \biggr | \\
   &\le& 
   \biggl | \| X (\hat \beta - \beta^0 ) \|_2 -m^0 \biggr | + 
    \biggl | m^0 - \sqrt n \| \Sigma_0^{1/2} ( \beta^* - \beta^0) \|_2 \biggr | \\
    &\le &  \biggl | \| X (\hat \beta - \beta^0 ) \|_2 -m^0 \biggr | +  (\gamma + \sqrt {2 \log n/n } ) \sqrt n \| \Sigma_0^{1/2} ( \beta^* - \beta^0) \|_2 + 4 \sqrt {\log 2} . 
    \end{eqnarray*}
    Apply Lemma \ref{concentration.lemma} again to finalize the result.
    \hfill $\sqcup \mkern -12mu \sqcap$
    
    \subsection{Proof of Theorem \ref{noiseless.theorem}}

To establish Theorem \ref{noiseless.theorem}, we first need to study the minimizer $b^*$ in (\ref{b*.equation}).
The minimization
$$ \min \biggl \{ \| X b \|_2^2  : \ \| b_S \|_1- \| b_{-S} \|_1 = 1 \biggr \} $$
has non-convex constraints. If we fix the signs within $S$ of a possible solution $b$, one can reformulate it as
a convex problem with convex constraints. This is done in Lemma \ref{convex.lemma}. 
We then show that $b_j^* \not= 0$ for all $j \in S$ in Lemma \ref{not0.lemma}. This is important because
given the signs within $S$ of a potential solution $b$, we want the restrictions on these signs to be
non-active so that the Lagrangian formulation is of a similar form as the KKT conditions 
(\ref{KKT-noiseless.equation}) for the noiseless Lasso.
This Lagrangian form is then given in Lemma \ref{z*.lemma} with  Lemma \ref{zS.lemma} serving as a preparation.
The Lagrangian form of Lemma \ref{z*.lemma} with $S=S_0$ in a sense resembles the KKT conditions
(\ref{KKT-noiseless.equation}) when the active coefficients in the vector $ \beta_S^0$ have appropriate signs
and $| \beta_j^0 |$ is for $j \in S_0$ large enough. This
allows one to find a solution $\beta^*$ of the KKT conditions (\ref{KKT-noiseless.equation}) with 
the prescribed prediction
error. 

\subsubsection{Non-sparseness within $S$}
Our first step is to ascertain that  a solution 
 $$ b^* \in \arg \min_{b \in \R^p} \biggl \{ \| X b \|_2: \ \| b_{S}  \|_1 - \| b_{-S} \|_1 =1 \biggr\} $$
 can be found by searching over (at most)
$2^{|S|}$ convex problems with convex constraints. This is done in the next lemma,
where we also show that the equality constraint $\| b_S \|_1- \| b_{-S} \|_1 =1$ can be
replaced by an inequality constraint $\| b_S \|_1- \| b_{-S} \|_1 \ge 1$. 

\begin{lemma}\label{convex.lemma}
We have
\begin{eqnarray*}
 &\  & \min  \biggl \{ \| X b \|_2^2 : \ \| b_S \|_1 - \| b_{-S} \|_1 = 1 \biggr \} \\
  &= & \min  \biggl \{\| X b \|_2^2 : \ \| b_S \|_1 - \| b_{-S} \|_1 \ge 1 \biggr \} \\
& = & \min_{z_S \in \{ \pm 1 \}^{|S|} } \min_{b} \biggl \{  \| X b \|_2^2   : \ z_j b_S  - \| b_{-S} \|_1 \ge 1, \ z_j b_j \ge 0 \ \forall \ j \in S  \biggr \} . 
\end{eqnarray*} 
\end{lemma}

{\bf Proof of Lemma \ref{convex.lemma}.}
To show that the equality constraint can be turned into an inequality constraint let us consider some $b \in \R^p$ for which it holds that $\| b_S \|_1- \| b_{-S} \|_1 =c $, where $c$ is a constant bigger than 1.
Let $\tilde b := b / c $. Then
$$\| \tilde b_S \|_1 - \| \tilde b_{-S} \|_1 = \biggl ( \| b_S \|_1- \| b_{-S} \|_1 \biggr ) / c = 1 . $$
Moreover
$$ \| X \tilde b \|_2 = \| X b \|_2 / c < \| X b \|_2 . $$
Thus the first equality of the lemma must be true.

We now show the second equality of the lemma. If for some $z_S \in \{ \pm 1 \} $ it holds that
 $ z_j b_j \ge 0$ for all $j \in S$, we have $z_S^T b_S= \| b_S \|_1$. Conversely,
 if we define for $j \in S$ with $b_j \not=0 $, $z_j := b_j / | b_j |$ as the sign of
$b_j$, and define $z_j \in \{ \pm 1 \} $ arbitrarily for $j \in S$ with $b_j = 0$, then we
have $z_j b_j \ge 0$ for all $j \in S$. Thus
$$
   \biggl \{ b: \ \| b_S \|_1 - \| b_{-S} \|_1 \ge 1 \biggr \} =   \cup_{z_S \in \{ \pm 1 \}^{|S|} } \biggl \{ b:  \ z_S^T b_S - \| b_{-S} \|_1 \ge 1 , \  z_j b_j \ge 0  
 \biggr \} .
 $$
\hfill $\sqcup \mkern -12mu \sqcap$

We establish in the next lemma that sign constraints on $b_S^*$ are not active: 
$b_S^*$ is so to speak maximally non-sparse. 
We assume that
$\hat \kappa^2 (S)>0$, so for $S=S_0$ we implicitly assume Condition \ref{inverse.condition}.

\begin{lemma}\label{not0.lemma}
Suppose that $\hat \kappa (S) \not= 0$. Then for any minimizer $b^*$
of the problem
$$ \min \biggl \{ \| X b \|_2 : \ \| b_S \|_1 - \| b_{-S} \|_1 = 1 \biggr \} $$
it holds that
$b_j^* \not= 0 $ for all $j \in S $.
\end{lemma}

\begin{remark}\label{hamza.remark}
A (very) special case  of Lemma \ref{not0.lemma} is the minimization problem
$$ {\rm b}_S^* \in  \argmin \biggl \{ \| b_S \|_2^2 : \ \| b_S \|_1= 1 \biggr \} .$$
Clearly the solution has $|{\rm b}_j^*|= 1/ |S| \not= 0 $ for all $j \in S$.
More generally, for the case without  ``$b_{-S}$-part" one can 
apply a geometric argument to show
that whenever $X_S^T X_S$ is non-singular
$${\rm b}_S^* \in   \argmin \{ \| X b_S \|_2: \ \| b_S \|_1 = 1 \} $$
must have all its components in $S$ nonzero.
Indeed , let $r:= \| X  {\rm b}_S^* \|_2$. Then $r >0$ by the   non-singularity of
$X_S^T X_S$.
Let
${\cal E}$ be the ellipsoid 
${\cal E} : = \{b_S  : \| X   b_S \|_2  \le r \}$ and ${\cal B} := \{ b_S : \| b_S \|_1 \le  1 \}$.
It is easy to see that ${\cal E}$ must be included in ${\cal B}$. Now
${\rm b}_S^*$  is a point on the boundary of both ${\cal E}$ and ${\cal B}$, so any supporting hyperplane to 
and ${\cal B}$ must also be supporting to ${\cal E}$. The key observation is that any point on the boundary of ${\cal E}$ 
has a unique supporting hyperplane (given by the gradient of the quadratic form); and that points on the boundary of 
${\cal B}$  that have a unique supporting hyperplane are exactly those points with no zero entry.
\end{remark}

{\bf Proof of Lemma \ref{not0.lemma}.} 
We use the representation of Lemma \ref{convex.lemma}.
Let $z_{S}^* \in \{ \pm 1 \}^{|S|}$ satisfy $z_{S}^{*T} b_S^* = \| b_S^* \|_1$ and
$z_j^* b_j^* \ge 0 $ for all $j \in S $. Then $b^*$ is a solution of the
convex minimization problem with (linear and) convex constraints
$$ \min \biggl \{ \| X b \|_2^2: \ z_S^{*T} b_S - \| b_{-S} \|_1 \ge 1 , z_j^* b_j \ge 0 , \ 
\forall \ j \in S \biggr \} .$$
Note that in the minimization, one may replace the inequality constraint
$z_S^{*T} b_S - \| b_{-S} \|_1 \ge 1 $ by an inequality constraint
$z_S^{*T} b_S - \| b_{-S} \|_1 = 1 $. This follows from the same arguments as used
in the proof of Lemma \ref{convex.lemma}. A reason to replace
the equality constraint by an inequality constraint is that
the restrictions become convex. 

The solution of the convex problem with convex constraints can be found
using
Lagrange multipliers $\tilde \lambda $ and $\mu_S$, where $\tilde \lambda \ge0$ and where $\mu_S$ is
an $|S|$-vector with non-negative entries.  
The Lagrangian formulation is
$$ \min \biggl \{ \| X b \|_2^2   +2 \tilde \lambda \biggl ( \| b_{-S} \|_1 - z_S^{*T} b_S -1 \biggr ) 
-2 \sum_{j \in S} \mu_{j,S} z_j^* b_j  \biggr \} . $$
Because the inequality constraint can be replaced by an equality
constraint, we know that in fact $\tilde \lambda >0$.
The Lagrangian formulation has has KKT conditions
$$ X^T X  b^* = \tilde \lambda z^* + {\rm diag} (\mu_S ) z_S^* , $$
where $z_{-S}^*$ is an element of the sub-differential
$$ - \partial \| b_{-S}^* \|_1 = \biggl \{ z_{-S}: \ \| z_{-S} \|_1 \le 1 , \ z_{-S}^T b_{-s}^* = 
- \| b_{-S}^* \|_1 \biggr \} . $$
It follows that for $j \in S$
$$ b_j^* \not= 0 \ \Rightarrow  \mu_{j,S} = 0 . $$
Let ${\cal N}:= \{j \in S :\  b_j^*= 0  \}$.
Then we have by the above argument
$$( X^T X b^*)_{-{\cal N}}=\tilde  \lambda z_{-{\cal N}}^* $$
$$( X^T X  b^* )_{{\cal N}} = \tilde \lambda z_{\cal N}^* + {\rm diag} ( \mu_{\cal N} ) z_{\cal N}^* . $$
The tangent plane of $\{ b: \ \| X b\|_2 = \| X b^{*}\|_2  \}$ at $b^*$ is
$${\cal U}:=  \{ u=b^* +v: \ v^T X^T X  b^* = 0  \} . $$
The idea of the proof is now to take an element $u= b^* + t v$ in this tangent plane with $t >0$ and
with $v_j \not= 0$ for at least one  $j \in {\cal N}$ and such that $v_j\not= 0 $ has the same sign as
$b_j^*$ for all $j \in S \backslash {\cal N} $. For $j \notin S$ we take $v_j=0$.
Then $\tilde b := b^* + tv$ has $\| \tilde b_S \|_1 - \| \tilde b_{-S} \|_1 > 1$
and this leads for a suitable scale $t$ to 
$${ \| X \tilde b \|_2\over 
\| \tilde b_S \|_1 - \| \tilde b_{-S} \|_1 } <\| X  b^* \|_2 . $$

Let us now work out this idea. 
It cannot be true that $b_j^*=0$ for all $j \in S$ as $\| b_S^* \|_1 \ge 1$. Hence $S \backslash {\cal N} \not=
\emptyset$. 
Take (for example) $v_j= z_j^* $ for all $j \in S \backslash {\cal N}  $.
Then 
$$v_{S \backslash {\cal N}}^T z_{S \backslash {\cal N}}^* =z_{S \backslash {\cal N}}^{*T} z_{S \backslash {\cal N}}^*=
| S \backslash {\cal N} | .$$
Now $\tilde \lambda >0$ and the entries of $\mu_{{\cal N}} $ are all positive as well (since $\mu_j =0$ for some
$j \in {\cal N}$ would imply $b_j^*=0$ for this $j$, which is not possible by the definition of ${\cal N } $). 
Therefore we can choose
$$v_{{\cal N} }^T ( \tilde \lambda z_{{\cal N}}^* + {\rm diag} (\mu_{{\cal N}}) z_{{\cal N}}^* ) =- \tilde \lambda | S \backslash {\cal N} | . $$
Then at least one entry of $v_{\cal N}$ has to be non-zero and moreover
\begin{eqnarray*}
 v^T X^T X  b^* &=& \tilde \lambda v_{S \backslash {\cal N}}^T z_{S \backslash {\cal N}}^* + v_{\cal N}^T ( \tilde \lambda z_{\cal N}^* + {\rm diag} (\mu_{\cal N}) z_{\cal N}^* ) \\
& = & \tilde \lambda | S \backslash {\cal N} |   -  \tilde \lambda | S \backslash {\cal N} |\\
&=&0 . 
\end{eqnarray*}
We thus have for all $t >0$
\begin{eqnarray*}
\| X (b^* + t v) \|_2^2 &=&
\| X  b^*\|_2^2  + t^2 \| X  v\|_2^2 . 
\end{eqnarray*}
Moreover
\begin{eqnarray*}
 \| b_S^* + t v _S\|_1  &= &\| b_{S \backslash {\cal N} }^* \|_1 + t \| v _{S \backslash {\cal N}} \|_1 +
t \| v_{\cal N} \|_1 \\
&=&  \| b_S^* \|_1 + t \| v \|_1. 
\end{eqnarray*}
Therefore
\begin{eqnarray*}
 \| b_S^* +  t v_S \|_1 - \| b_{-S}^* \|_1 &=&  \| b_S^* \|_1 -
\| b_{-S}^* \|_1 +  t \| v \|_1\\
&= &1+   t \| v \|_1 . 
\end{eqnarray*}
It follows that
\begin{eqnarray*}
& \ & { \| X  (b^* + t v) \|_2^2 \over
(\| b_S^*+ t v_S  \|_1 - \| b_{-S}^* \|_1 )^2 }\\
& =&  { \| X  b^*\|_2^2  + t^2 \| X  v \|_2^2 
\over ( 1+   t \| v \|_1)^2 } .
\end{eqnarray*}
Define
\begin{eqnarray*}
 {\rm A}& :=  &
 \| X  b^* \|_2^2 + t^2 \| X  v\|_2^2 -
\| X b^* \|_2^2  ( 1+   t \| v \|_1)^2 \\
&=& t^2 \| X v \|_2^2 - 2t \| X  b^* \|_2^2  \| v \|_1 -
 t^2 \| X  b^*\|_2^2  \| v \|_1^2  \\
& = & t^2 (\| X v \|_2^2 - \| X  b^* \|_2^2  \| v \|_1^2 ) 
- 2t \| X b^* \|_2^2   \| v\|_1^2 .
\end{eqnarray*}
We will show that for suitable $t>0$ the constant ${\rm A}$ is strictly negative.
This means
$$ \| X ( b^*+ t v) \|_2^2    <
\| Xb^*\|_2^2  (  \|b_S^* +  t v_S \|_1- \| b_{-S}^* \|_1 )^2 $$
and so we arrive at a contradiction.
To show ${\rm A}<0$ we distinguish two cases. If
$$\| X  v \|_2^2 \le  \| X b^* \|_2^2  \| v\|_1^2 $$
then ${\rm A} <0$ for all $t>0$. If
$$\| X  v \|_2^2 >  \| X b^* \|_2^2  \| v\|_1^2 $$
then ${\rm A} <0$ for all $t$ satisfying
$$0<t < { 2 \| X b^* \|_2^2   \| v \|_1^2  \over 
\| X  v\|_2^2  - \| X b^* \|_2^2 \| v \|_1^2 }   . $$
Here we used the assumption that $\| X b^* \|_2^2 >0$ so that the above
right hand side is indeed strictly positive. 
\hfill $\sqcup \mkern -12mu \sqcap$

\subsubsection{Lagrangian form} 
We now present the Lagrangian form given the signs within the set $S$ and given that within
the set $S$ the solution has non-zero entries.
Let for each $z_S \in \{ \pm 1 \}^{|S|}$
$$
 b^* (z_S)\in  \arg \min  \biggl \{ \| X b \|_2^2  : \ z_S^T b_S - \| b_{-S} \|_1 \ge 1, z_j b_j \ge 0, \ \forall \  j \in S   \biggr \} . 
$$
Define
$${\cal Z}_S := \biggl \{ z_S \in \{ -1 , 1 \}^{|S| }: \ z_j b_j^* (z_S) >0 \ \forall \ j \in S \biggr \} . $$

\begin{lemma} \label{zS.lemma} We have for all $z_S \in {\cal Z}_S $
$$ X^T X  b^* (z_S) = { z^* (z_S) \| X b^{*} (z_S)\|_2^2  } $$
where $z_S^* (z_S) = z_S$ and 
$z_{-S}^*(z_S) \in - \partial \| b_{-S} ^* (z_S) \|_1 $. 
\end{lemma}

{\bf Proof of Lemma \ref{zS.lemma}.}  To prove this result it is useful to repeat some arguments
of the proof of Lemma \ref{not0.lemma}. 
The
convex minimization problem with (linear and) convex constraints
$$
\min  \biggl \{ \| X b \|_2^2  : \ z_S^T b_S - \| b_{-S} \|_1 \ge 1, z_j b_j \ge 0, \ \forall \  j \in S   \biggr \} 
$$
can be solved using Lagrange multipliers $\tilde \lambda $ and $\mu_S$, where
$\tilde \lambda>0$ and $\mu_S$ is an $|S|$-vector with non-negative entries.
The Lagrangian formulation is 
$$ \min \biggl \{ \| X b \|_2^2  + 2 \tilde \lambda \biggl ( \| b_{-S} \|_1 - z_S^T b_S -1 \biggr ) - 2\sum_{j \in S}
\mu_{j,S} z_j b_j   \biggr \} .$$
This has KKT conditions
$$ X^T X  b^* (z_S) = \tilde \lambda z^* + {\rm diag} (\mu_S) z_S ,$$
where $z_S^* = z_S$ and $z_{-S}^* = z_{-S}^* (z_S) $ depends on $z_S$ and is an element
of the sub-differential 
$$- \partial \| b_{-S} ^* (z_S) \|_1= \biggl \{ z_{-S}: \ \| z_{-S} \|_{\infty} \le 1 , \ 
z_{-S}^T b_{-S}^* (z_S )= - \| b_{-S}^* \|_1  \biggr \}  . $$
It follows that for $j \in S$
$$b_j^* (z_S) \not= 0  \Rightarrow \mu_{j,S} =0. $$
The assumption that $z_S \in {\cal Z}_S$ thus gives $\mu_S=0$. The KKT conditions
then read
$$ X^T X  b^* (z_S) = \tilde \lambda z^* . $$
One sees that 
$$ 1= z^{*T} b^* (z_S) =  b^{*T} (z_S) X^T X  b^* (z_S)/\tilde \lambda = \| X b^* (z_S) \|_2^2 / \tilde \lambda  . $$
This gives
$$ \tilde \lambda = {  \| X  b^* (z_S) \|_2^2 } \ . $$
\hfill $\sqcup \mkern -12mu \sqcap$

We apply the above lemma with $z_S:=  \partial \| b_S^*\|_1 $. This gives the following result.

\begin{lemma}  \label{z*.lemma} Suppose $\hat \kappa (S) \not= 0$. Let
$$ b^* \in \arg \min \biggl \{ \| X b \|_2^2  : \ \| b_S \|_1 - \| b_{-S} \|_1 =1 \biggr \}  $$
Then
$$X^T X  b^*  = z^*  \| X b^* \|_2^2  . $$
where
$z_S^* = \partial \| b_S^* \|_1$ and $z_{-S}^* \in - \partial \| b_{-S} ^*  \|_1$. 
\end{lemma}

{\bf Proof of Lemma \ref{z*.lemma}.} 
By Lemma \ref{not0.lemma},
for each 
$$ b^* \in \arg \min \biggl \{ \| X b \|_2^2  : \ \| b_S \|_1 - \| b_{-S} \|_1 = 1 \biggr \} $$
it holds that $b_j^* \not= 0 $ for all $j \in S$. We can therefore define
$z_j^*:= b_j^* / | b_j^* | $ for all $j \in S$ and then $z_S^* = \partial \| b_S^* \|_1  \in {\cal Z}_S$. 
The result now follows from Lemma \ref{zS.lemma}.
\hfill $\sqcup \mkern -12mu \sqcap$ 

\subsubsection{Finalizing the proof of Theorem \ref{noiseless.theorem}.} 

With the help of Lemma \ref{z*.lemma} we are now in the position to prove  Theorem \ref{noiseless.theorem}. 

{\bf Proof of Theorem \ref{noiseless.theorem}.} 
Let $b^*$ and $z^*$ be as in Lemma \ref{z*.lemma}, with $S= S_0$.
Define
$$ \beta^{\prime} = \beta^0  -  { b^* s_0 \over  \hat \kappa^2 (S_0)} 
{ \lambda^* \over n} . $$ 
Then
\begin{eqnarray*}
X^T X  ( \beta^{\prime} - \beta^0) & = & - { \lambda^* X^T X  b^* s_0  \over n  \hat \kappa^2 (S_0) } \\
 & = & -{ \lambda^* X^T X  b^* \over \| X b^* \|_2^2}  \\
 & = &  - \lambda^* z^* .
 \end{eqnarray*} 
Let $S_* := \{ j : \ b_j^* \not= 0 \} $. Then by Lemma \ref{not0.lemma}, $S_0 \subset S_*$.
Furthermore
$$ z_j^* \beta_j^{\prime} = \begin{cases}  z_j^* \beta_j^0 - \lambda z_j^* b_j^* / \| X b^* \|_2^2  >0 & j \in S_0 \cr  
- \lambda^* z_j^*  b_j^* / \| X b^* \|_2^2  >0 & j \in S_* \backslash S_0 \cr 
0 & j \notin S_* \cr  \end{cases} . $$
It follows that $z^* \in \partial \| \beta^{\prime} \|_1 $.
Thus, $\beta^{\prime}=: \beta^*$ is a solution of the KKT conditions 
(\ref{KKT-noiseless.equation}) with $\zeta^*= z^* $.
It holds moreover that
$$ \| X ( \beta^* - \beta^0) \|_2^2  = {  \lambda^{*2}  \| X b^{*} \|_2^2  \over\| X b^* \|_2^4 } =
{ \lambda^{*2} s_0  \over n \hat \kappa^2 (S_0) }.$$
\hfill $\sqcup \mkern -12mu \sqcap$

\subsection {Proof of Theorem \ref{TV-compatibility.theorem}.} 
The proof of Theorem \ref{TV-compatibility.theorem} consists of several steps.
First we note that, given the sizes of its jumps, the total variation of a function is the smallest when
this function is decreasing or increasing. This is stated in Lemma
\ref{trivial1.lemma} as a trivial fact. As a consequence, if one subtracts from
an arbitrary function value - or minus this value -  the total variation, the result will be at most
the average of the absolute values. This is shown in Lemma \ref{trivial2.lemma}. 
Lemma \ref{trivial2.lemma} is then applied at each jump separately, as $\| b_S\|_1 - \| b_{-S}\|_1$
in this example amounts to subtracting at each jump some total variation to the left or
to the right of this jump. Lemma \ref{trivial3.lemma} shows how this works for one jump.
Then Theorem \ref{TV-compatibility.theorem} is in part proved by applying
this lemma to each jump. This leads to a lower bound for $\hat \kappa^2 (S)$.
The proof is completed by showing that this lower bound is achieved by the vector $b^*$ as given
in Theorem \ref{TV-compatibility.theorem}. 

  For $f \in \R^n$ we define the ordered vector
  $$f_{(1)} \le \cdots \le f_{(n)} , $$
  with arbitrary ordering within ties.
  
  \begin{lemma} \label{trivial1.lemma} It holds that
  $$ TV(f) \ge f_{(n)}- f_{(1)} $$
  with equality if $f$ is increasing or decreasing.
    \end{lemma}
    
    {\bf Proof of Lemma \ref{trivial1.lemma}.} Trivial.
    \hfill $\sqcup \mkern -12mu \sqcap$
      
  \begin{lemma} \label {trivial2.lemma} It holds for any $j \in \{1 , \ldots , n\} $ that
  $$ f_j - {\rm TV} (f) \le f_{(1)} \le {1 \over n} \sum_{i=1}^n |f_i| , $$
  and
  $$ - f_j - {\rm TV} (f) \le - f_{(n)} \le {1 \over n} \sum_{i=1}^n |f_i| . $$
   \end{lemma} 
   
   {\bf Proof of Lemma \ref{trivial2.lemma}.} 
   We have from Lemma \ref{trivial1.lemma} that ${\rm TV} (f) \ge f_{(n)} - f_{(1)} $.
   Moreover, $f_j \le f_{(n)} $. Thus
   \begin{eqnarray*}
  f_j - {\rm TV} (f) &\le & f_j - (f_{(n)} - f_{(1)}) \\
  & \le & f_{(n)} - (f_{(n)} - f_{(1)} )\\
  &=&  f_{(1)} .
  \end{eqnarray*} 
  \fbox{Case 1:} if $f_{(1)} < 0 $ obviously $f_{(1)} < {1 \over n} \sum_{i=1}^n |f_i|$.\\
  \fbox{Case 2:} if $f_{(1)} \ge 0$ then $f_i \ge 0 $ for all $i$ and then
  $$f_{(1)} \le \sum_{i=1}^n f_i / n = \sum_{i=1}^n |f_i|  / n . $$
  In the same way
    \begin{eqnarray*}
  -f_j - {\rm TV} (f) &\le & -f_j - (f_{(n)} - f_{(1)}) \\
  & \le & -f_{(1)} - (f_{(n)} - f_{(1)} )\\
  &=&  -f_{(n)} .
  \end{eqnarray*} 
 \fbox{Case 1:} if $f_{(n)} > 0 $ then $-f_{(n)} < {1 \over n} \sum_{i=1}^n |f_i|$.\\
 \fbox{Case 2:} if $f_{(n)} \le 0$ then $f_i \le 0 $ for all $i$ and then
  $$-f_{(n)} \le - \sum_{i=1}^n f_i / n = \sum_{i=1}^n |f_i|  / n . $$
  \hfill $\sqcup \mkern -12mu \sqcap$
    
  \begin{lemma}\label{trivial3.lemma} Let $f \in \R^n$ with total variation
  ${\rm TV} (f) = \sum_{i=2}^n | f_i - f_{i-1} |$ and $g \in \R^m$ with total variation
  ${\rm TV} (g) = \sum_{i=2}^m | g_i - g_{i-1} |$. Then for any $j \in \{1 , \ldots , n\}$ and
  $k \in \{1 , \ldots, m \}$
  $$ | f_j - g_k | - {\rm TV} (f) - {\rm TV} (g) \le {1 \over n} \sum_{i=1}^n | f_i | +
  {1 \over m} \sum_{i=1}^m | g_i | . $$
   \end{lemma}
   
  {\bf Proof of Lemma \ref{trivial3.lemma}.} 
  Suppose without loss of generality that $f_j \ge  g_k$. Then by Lemma \ref{trivial2.lemma}
  \begin{eqnarray*}
   | f_j - g_k | - {\rm TV} (f) - {\rm TV} (g)&=& \underbrace{(f_j - {\rm TV} (f))}_{\le \sum_{i=1}^n | f_i | /n }  +
   \underbrace{(- g_k - {\rm TV} (g))}_{\le \sum_{i=1}^m | g_i | / m }  \\
   & \le & {1 \over n} \sum_{i=1}^n | f_i | +
  {1 \over m} \sum_{i=1}^m | g_i | .
    \end{eqnarray*}
    \hfill $\sqcup \mkern -12mu \sqcap$
    
  {\bf Proof of Theorem \ref{TV-compatibility.theorem}.}
  Let for $j=2 , \ldots , s$, $u_j \in \Nat$ satisfy $1\le u_j \le d_j -1 $.
  We may write for $f = X b$,
  \begin{eqnarray*}
  & \ & 
   \| b_S \|_1- \| b_{-( S \cup \{ 1 \} ) } \|_1 \\
   &=&   | f_{d_1+1} - f_{d_1} | - \sum_{i=2}^{d_1} | f_i - f_{i-1} |
 - \sum_{i=d_1+2}^{d_1+u_2} 
  | f_i- f_{i-1} | \\ &+& | f_{d_1+d_2+1} - f_{d_1+d_2} |-
   \sum_{i=d_1+u_2+1}^{d_1+d_2} 
  | f_i- f_{i-1} | - \sum_{i=d_1+d_2 +2}^{d_1+d_2+u_3} |f_i - f_{i-1} |  \\
  & \ & \cdots \\ 
  &+&  | f_{d_1+\cdots + d_{s-1}+1} - f_{d_1+\cdots + t_{s-1} } |\\ &-&
   \sum_{i=d_1+\cdots + d_{s-2}+ u_{s-1}+1}^{d_1+ \cdots + d_{s-1}} 
  | f_i- f_{i-1} | 
- \sum_{i=d_1+\cdots + d_{s-1} +2}^{d_1 + \cdots + d_{s-1} + u_s }   | f_i- f_{i-1} | \\
  &+&  | f_{d_1+\cdots + d_s+1} - f_{d_1+\cdots + d_s } |\\ &-&
   \sum_{i=d_1+\cdots + d_{s-1}+ u_s+1}^{d_1+ \cdots + d_s} 
  | f_i- f_{i-1} | 
- \sum_{i=d_1+\cdots + d_{s} +2}^{n}   | f_i- f_{i-1} | 
   \end{eqnarray*}
   \begin{eqnarray*}
   & \le& { 1 \over d_1} \sum_{i=1}^{d_1} | f_i | + { 1 \over u_2} \sum_{i=d_1+1}^{d_1+u_2 }  | f_i  | \\
   &+& { 1 \over d_2 - u_2} \sum_{i=d_1+u_2+1}^{d_1+ d_2}  | f_i| +{1 \over u_3}
   \sum_{i=d_1+ d_2+1 }^{d_1 + d_2 + u_3} |f_i| \\
   & \ & \cdots \\
   & +& {1 \over d_{s-1} - u_{s-1}} \sum_{i=d_1 + \cdots + d_{s-2} + u_{s-1}+1}^{d_1 + \cdots + d_{s-1}} | f_i  | 
   + {1 \over u_{s}} \sum_{i=d_1 + \cdots + d_{s-1}+1}^{d_1 + \cdots + d_{s-1} + u_s  }   | f_i |\\
   &+& { 1 \over d_s -u_s} \sum_{i=d_1 + \cdots + d_{s-1} + u_{s}+1}^{d_1 + \cdots + d_{s}} | f_i  | +
   { 1 \over d_{s+1}  } \sum_{i=d_1 + \cdots + d_{s}+1}^{n}   | f_i |
   \end{eqnarray*}
   \begin{eqnarray*}
   &\le& \sqrt { {1 \over d_1} + {1 \over u_2} + {1 \over d_2 - u_2 } + \cdots + {1 \over d_{s-1}- u_{s-1}}+
   {1 \over u_{s}}  + {1 \over d_s - u_s } + {1 \over d_{s+1} }   } \\
   & \ & \times \sqrt { \sum_{i=1}^n | f_i |^2 } ,
    \end{eqnarray*}
    where in the first inequality we applied
    Lemma \ref{trivial3.lemma} and the second one follows from the Cauchy-Schwarz inequality. 
The assumption that for all $j \in \{ 2, \ldots , s\}$ that $d_j $ is even allows us to take $u_j = d_j / 2$
to arrive at 
$$ \kappa^2 (S) \ge 
{ s+1 \over { n \over d_1}  + \sum_{j=2}^s { 4 n\over d_j } + { n \over d_{s+1} }} . $$
     
     Now for the reverse inequality, let $\tilde b$ be given as in the theorem and
     and $\tilde f:= X \tilde b $. Then $\tilde f$ is equal to
  $${\tilde f}_i=\begin{cases} -{n \over d_1} & i=1 , \ldots , d_1 \cr 
   {2 n\over d_2} & i=d_1+1 , \ldots , d_1 + d_2 \cr  
   \vdots & \ \cr (-1)^s {2n \over d_s} & i= \sum_{j=1}^{s-1} d_j +1 , \ldots ,\sum_{j=1}^s d_j  \cr 
   (-1)^{s+1} {n \over d_{s+1} } & i= \sum_{j=1}^s d_j +1 , \ldots , n \cr \end{cases} . $$
By the definition of $\tilde f= X \tilde b$, 
      \begin{eqnarray*}
      \| \tilde b_{S} \|_1= 
    \sum_{j=1}^s | {\tilde f}_{d_j+1} - {\tilde f}_{d_j} | &=& { n \over d_1} + {2 n\over d_2} \\
    &+&  {2n \over d_2} + {2n \over d_3} \\
    & \ & \vdots \\
    &+& {2n \over d_{s-1}} + {2n \over d_s} \\
    &+& { 2n \over d_s} + { n \over d_{s+1} } \\
    &=& { n \over d_1} +  \sum_{j=2}^s { 4n\over d_j} + {n  \over d_{s+1} } ,
    \end{eqnarray*}
    and also
    \begin{eqnarray*}
    \sum_{i=1}^n {\tilde f}_i^2 &=&
    { d_1 } {\tilde f}_{t_1}^2 + \cdots + { d_{s+1} } {\tilde f}_{d_{s+1}}^2 \\
    &=& { n^2 \over d_1} + 4 \sum_{j=2}^s { n^2 \over d_j} + {n^2 \over d_{s+1} } . 
    \end{eqnarray*}
    Note also that
    \begin{eqnarray*}
  & \ &  \| \tilde b_{-(S \cup \{1 \} ) } \|_1\\
   &=& \sum_{i=2}^{d_1} | {\tilde f}_{i}- {\tilde f}_{i-1} | +  \sum_{i=d_1+2}^{d_2} |{\tilde f}_{i} - {\tilde f}_{i-1} | + \cdots + \sum_{i=d_1 + \cdots + d_s +2}^{n} | {\tilde f}_i - {\tilde f}_{i-1}|  \\
   &=&0 
   \end{eqnarray*} 
    It follows that
    \begin{eqnarray*}
     { (s+1) \| X \tilde  b \|_2^2/n  \over  ( \| \tilde b_S \|_1 - \| \tilde b_{-(S \cup \{1 \} ) } \|_1 )^2 } & = & 
    {  \sum_{i=1}^n {\tilde f}_i^2/n  \over 
      \biggl ( \sum_{j=1}^s | {\tilde f}_{d_j+1} - {\tilde f}_{d_j} | \biggr )^2 } \\
     &=& { s+1 \over { n \over d_1}  + \sum_{j=2}^s { 4n \over d_j } + { n \over d_{s+1} }} . 
     \end{eqnarray*} 
     \hfill $\sqcup \mkern-12mu \sqcap$

   \subsection{Proof of Theorem \ref{noisy.theorem}.} 
    
      To prove Theorem \ref{noisy.theorem}, we first establish the Lagrangian form of the minimization
      problem where we have the convex constraint $ z_{S_0}^{*T}(\bar v)  b_{S_0} -
      \| W b_{-S_0} \|_1 \ge 1 $. Then we recall the projections and we introduce a subset ${\cal T}$
      of the underlying probability space where the lower bound of Theorem \ref{noisy.theorem}
      holds. The latter is shown in Lemma \ref{deterministic2.lemma}.
      Finally, we show that the subset ${\cal T}$ has large probability.

      \subsubsection{Lagrangian form}
       Recall for $w  \in {\cal W} (\bar v) $ the convex problem with linear and convex constraints
   $$b(w ) \in \argmin \biggl \{ \| X b \|_2^2: \ z_{S_0}^{*T}  (\bar v) b_{S_0} - \| W b_{-S_0} \|_1 \ge 1 \biggr \} . $$
      Note that here we do not require the positivity constraint
      $ z_j^{*T}(\bar v)  b_j \ge 0 $ for all $j \in S_0$.  The next lemma gives its Lagrangian form.
      This  form plays in the proof of Theorem \ref{noisy.theorem} 
      the same role as in the proof of Theorem \ref{noiseless.theorem} for the noiseless version.
      We also show that for $w \in {\cal W} (\bar v) $ the minimum $\| X b(w) \|_2^2 $ is not larger than $\| X  b^* (\bar v)  \|_2^2 $
      (recall that by definition $\hat \kappa^2 ( 1+ \bar v, S_0) = s_0 \| X  b^* (\bar v)  \|_2^2 /n $). 
   
   \begin{lemma} \label{linear.lemma}
   We have
   $$ X^T X  b(w)= \| X b(w) \|_2^2 W z(w) , $$
    with 
    $$z_{S_0} (w) =  z_{S_0}^* (\bar v)  , \  z_{-S_0} (w) \in -\partial \| b_{-S_0} (w) \|_1 .$$
    Moreover, for $w \in {\cal W} (\bar v)$
    $$s_0 \| X b(w) \|_2^2 /n\le \hat \kappa^2 ( 1+ \bar v , S_0)  . $$
    \end{lemma}
    
    {\bf Proof of Lemma \ref{linear.lemma}.}
   The problem 
   $$\min \biggl \{ \| X b \|_2^2: \ z_{S_0}^{*T}(\bar v)   b_{S_0} - \| W b_{-S_0} \|_1 \ge 1 \biggr \}  $$
   has Lagrangian
   $$ X^T X   b(w) = \tilde \lambda  W z (w) $$
   with $z_{S_0} (w) = z_{S_0}^* (\bar v)  $ and $z_{-S_0} (w) \in -\partial \| b_{-S_0} (w) \|_1 $. 
   Moreover
   $$ \| X b(w) \|_2^2 = \tilde \lambda b(w)^T W z(w) =  z_{S_0}^{*T} (\bar v)  b_{S_0} - \| W b_{-S_0} \|_1 =1 $$
   because the minimum is reached at the boundary. 
 So
   $$ \tilde \lambda = \| X b(w) \|_2^2 . $$ 
   To obtain the second statement of the lemma, 
   we use similar arguments as in the proof of Lemma \ref{improve.lemma}. We have
   $$ \| X b(w) \|_2 = \min_{b \in \R^p} \biggl \{ 
   { \| X b \|_2 \over  z_{S_0}^{*T} (\bar v) b_{S_0}- \| W_{-S_0} b_{-S_0} \|_1 } : \ 
   z_{S_0}^{*T}(\bar v) b_{S_0} - \| W_{-S_0} b_{-S_0} \|_1 > 0 \biggr \} $$
   But for $w \in {\cal W}$ and $\bar w := 1 + \bar v $, we know 
   $$  \| W b_{-S_0} \|_1 \le \| \bar W  b_{-S_0} \|_1 $$
   and so
   $$ z_{S_0}^{*T} (\bar v) b_{S_0} - \| W b_{-S_0} \|_1 > 
    z_{S_0}^{*T}(\bar v)  b_{S_0 } -  \| \bar W  b_{-S_0} \|_1 . $$
    Let
    $$A:= \biggl \{ b: \ z_{S_0}^{*T} (\bar v) b_{S_0} - \| W b_{-S_0} \|_1 >0 \biggr \}$$
    and
    $$B := \biggl \{ b:\  z_{S_0}^{*T} (\bar v) b_{S_0 } -  \| \bar W   b_{-S_0} \|_1>0 \biggr \} . $$
    Then $B \subset A $. 
    Hence
    \begin{eqnarray*}
    \| X b(w) \|_2  &=&
   \min_{b \in A}  { \| X b \|_2 \over z_{S_0}^{*T}(\bar v)  b_{S_0}- \| W b_{-S_0} \|_1 } \\
   & \le & \min_{b \in B}  { \| X b \|_2 \over z_{S_0}^{*T}(\bar v)  b_{S_0}- \| W b_{-S_0} \|_1 }\\
   & \le & \min_{b \in B} { \| X b \|_2 \over  z_{S_0}^{*T} (\bar v) b_{S_0}- \| \bar W b_{-S_0} \|_1 }\\
   &=& \| X b^* (\bar v) \|_2 \\
   &=& {\sqrt n \hat \kappa (1+ \bar v, S_0) \over \sqrt {s_0}  } . 
     \end{eqnarray*}
     \hfill $\sqcup \mkern -12mu \sqcap$

           \subsubsection{Projections}\label{projections2a.section}
           Recall the notation of Subsection \ref{projections2.section} and that moreover
      the diagonal elements of the matrix
   $(X_{S_0}^T X_{S_0} )^{-1} $ are denoted by $\{ u_j^2 \}_{j \in S_0} $.
   We write
   $$ \hat u_{S_0} := (X_{S_0}^T X_{S_0} )^{-1} X_{S_0}^T \epsilon  .$$
   We denote the projection of $\epsilon $ on the space spanned by the columns of $X_{S_0}$ by
    $$ \epsilon {\rm P} X_{S_0} := X_{S_0} (X_{S_0}^T X_{S_0} )^{-1} X_{S_0}^T \epsilon =
    X_{S_0} \hat u_{S_0}   $$
    and write
    $${\bf U} (S_0):= \| \epsilon {\rm P } X_{S_0} \|_2  .$$

     \subsubsection{Choice of $\lambda$} \label{lambda2a.section}
 Recall that we require that for some $t>0$    
  $$
     \lambda>  \| v_{-S_0} \|_{\infty} 
   \sqrt {2  (\log (2p) +t) } .
  $$

  \subsubsection{The set ${\cal T}$}\label{setT.section}
  Recall
   \begin{equation} \label{uv.equation}
   \bar u_j := u_j \sqrt {2  (\log (2p) +t) }/ \lambda ,\ j \in S_0 ,  \ \bar v_j := v_j \sqrt {2  (\log (2p) +t) }/ \lambda , \
   j \notin S_0 .
   \end{equation}
  Let
   ${\cal T}$ be the set
   \begin{eqnarray*}
   {\cal T}
    &:=&  \biggl \{ | \hat u_j | \le \lambda \bar u_j \ \forall j \in S_0 \biggr \} \\
   & \cap &   \biggl \{ |\hat v_j | \le \lambda \bar v_j  \ \forall j \notin S_0 \biggr \}
    \cap  \biggl  \{  {\bf U} (S_0)  \le \sqrt {  s_0 } + \sqrt {2x} \biggr \} . 
   \end{eqnarray*}
   We show in Subsection \ref{probabilistic2.section}
   that $\PP ( {\cal T}) \ge 1-  \exp [-t] - \exp[-x]$.

 \subsubsection{Deterministic part} 

   The idea is now to incorporate the noisy part of the KKT
     conditions for the noisy Lasso into a weighted sub-differential, creating in that way
     KKT conditions of the same for as the noiseless KKT conditions (see
     (\ref{fakeKKK.equation}) in the proof). To do so, we
     first put part of the noise in the vector $\beta^0$ without adding additional non-zeros. 
     This makes it possible not to change the sub-differential at $S_0$. 
     The rewriting of the KKT conditions make them resemble the Lagrangian form of
     Lemma \ref{linear.lemma}.
     
     We will use the
     KKT conditions (\ref{hatKKT.equation}) for $\hat \beta$: 
   $$ - X^T (Y-X \hat \beta)  = -\lambda \hat \zeta  , \ \hat \zeta \in \partial \| \hat \beta \|_1 .  $$

    \begin{lemma}\label{deterministic2.lemma}
     Suppose we are on the set ${\cal T}$ defined in Subsection
     \ref{setT.section}. Then under the conditions
     of Theorem \ref{noisy.theorem}
     $$ \| X (\hat \beta - \beta^0 )\|_n \ge {\lambda \sqrt {s_0} \over
     \sqrt n \hat \kappa ( 1+ \bar v , S ) } + \sqrt {2x}  $$
     \end{lemma}

       {\bf Proof of Lemma \ref{deterministic2.lemma}.}    
       Set
   $$ \hat \beta_{S_0}^0 := \beta^0 + \hat u_{S_0} , \ \hat \beta_{-S_0}^0 := 0 . $$
  Then
   \begin{eqnarray*}
    Y&= &X \beta^0 + \epsilon \\
    &=&  X_{S_0} \beta_{S_0}^0+ X_{S_0} \hat u_{S_0} + \epsilon {\rm A} X_{S_0}\\
    &=& X \hat \beta^0 + \epsilon {\rm A} X_{S_0} . 
    \end{eqnarray*}
    The KKT conditions (\ref{hatKKT.equation}) are
   $$ - X^T (Y-X \hat \beta)  = -\lambda \hat \zeta  . $$
   We have
   $$ Y- X \hat \beta= - X (\hat \beta - \hat \beta^0 )- \epsilon {\rm A} X_{S_0} . $$
   Therefore
  $$ -X^T ( Y - X \hat \beta) =  X^T X (\hat \beta - \hat \beta^0 ) -
  X^T ( \epsilon {\rm A} X_{S_0} ). $$
  But
  $$ X_{S_0}^T (\epsilon {\rm A} X_{S_0} ) =0 , $$
  and
  \begin{eqnarray*}
   X_{-S_0}^T ( \epsilon {\rm A} X_{S_0} )& = &X_{-S_0}^T- X_{-S_0}^T X_{S_0} 
  ( X_{S_0}^T X_{S_0} )^{-1} X_{S_0}^T \epsilon \\
  & =& (X_{-S_0} {\rm A} X_{S_0} )^T \epsilon . 
  \end{eqnarray*} 
  Hence the KKT conditions read
  $$ X^T X  ( \hat \beta - \hat \beta^0 ) = - \lambda \hat \zeta + \hat v , $$
  where
  $$ \hat v_{S_0} =0 , \ \hat v_{-S_0} = (X_{-S_0} {\rm A} X_{S_0} )^T \epsilon . $$
  Set $\hat S:= \{ j: \ \hat \beta_j \not= 0 \} $ and define  for all $j \in \hat S \backslash S_0$ 
   $$\hat w_j := 1  + \hat v_j / ( \lambda  \hat \zeta_j)   . $$
  By assumption (since we are on ${\cal T}$) $ |\hat v_j |< \lambda \bar v_j$.
  so $\hat w_j \ge 1- \bar v_j$ for all $j  \in \hat S \backslash S_0$.
  For $j \notin \hat S \cup S_0$  we define
  $$ \hat w_j := \max \{ | 1 + \hat v_j /  \lambda  | , 1- \bar v_j \} .$$
  Then for $j \notin \hat S \cup S_0$
  \begin{eqnarray*}
   \lambda \hat \zeta_j + \hat v_j &=&
 \lambda  |  \hat \zeta_j + \hat v_j / \lambda | {\rm sign} ( \hat \zeta_j + \hat v_j / \lambda ) \\
 &=& 
 \begin{cases} \hat w_j  {\rm sign} ( \hat \zeta_j + \hat v_j / \lambda )  , &  | \hat \zeta_j + \hat v_j / \lambda | \ge
 1- \bar v_j  
  \cr \hat w_j {  |  \hat \zeta_j + \hat v_j / \lambda | \over 1- \bar v_j } {\rm sign} ( \hat \zeta_j + \hat v_j / \lambda )& 
   | \hat \zeta_j + \hat v_j / \lambda |  \le 1- \bar v_j \cr  \end{cases}\\
   &=& \hat  w_j \tilde \zeta_j , 
   \end{eqnarray*} 
   where
  $$ \tilde \zeta_j :=  \begin{cases}  {\rm sign} ( \hat \zeta_j + \hat v_j / \lambda )  , &  | \hat \zeta_j + \hat v_j / \lambda | \ge
 1- \bar v_j  
  \cr  {  |  \hat \zeta_j + \hat v_j / \lambda | \over 1- \bar v_j } {\rm sign} ( \hat \zeta_j + \hat v_j / \lambda )& 
   | \hat \zeta_j + \hat v_j / \lambda |  \le 1- \bar v_j \cr  \end{cases} . $$
   
 One readily verifies that (on ${\cal T}$) $\hat w_j \le 1+ \bar v_j $ for all $j \notin S_0$. 
   Taking $\tilde \zeta_j = \hat \zeta_j$ for $j \in S\cup S_0$  we arrive at the KKT conditions
   \begin{equation}\label{fakeKKK.equation}
  X^T X ( \hat \beta - \hat \beta^0  ) = -
  \lambda   \hat {W}  \tilde  \zeta  , \ \tilde \zeta \in \partial \| \hat \beta \|_1 
  \end{equation}
  and where $\hat W = {\rm diag} (\hat w)$ with $\hat w \in {\cal W} (\bar v  )  $. 
  Let now $S_0^+ := \{ j \in S_0:\  z_j^*(\bar v)  b_j (\hat w) > 0 \} $ and
  $S_0^- := \{ j \in S_0:\  z_j^* (\bar v) b_j (\hat w) \le 0 \} $.
 Take
  $$ \beta^{\prime} = \hat \beta^0 - \lambda b_j (\hat w ) / \| X b (\hat w) \|_2^2 . $$
  \fbox{Case 1} Let $j \in S_0$.
   By our condition on $\beta^0$ we know that for $j \in S_0$, $|\beta_j^0 | >
  \lambda |  b_j (\hat w )|  / \| X b (\hat w) \|_2^2 + |\hat u_{S_0} |$,
  so
  $| \hat \beta_j^0 | \ge  |\beta_j^0 | - | \hat u_{S_0} | > \lambda |  b_j (\hat w )|  / \| X b (\hat w) \|_2^2 $.
  If $z_j^* (\bar v) =1 $ and $b_j (\hat w) >0 $, then $\hat \beta_j^0>0$ and 
  $$\beta_j^{\prime} =  | \hat \beta^0|  -  \lambda | b_j (\hat w )| / \| X b (\hat w) \|_2^2 > 0 . $$
  If $z_j^* (\bar v) =1 $ and $b_j (\hat w) \le 0$,  then $\hat \beta_j^0>0$ and  we have
  $$\beta_j^{\prime} = | \hat \beta_j^0| +  \lambda |b_j (\hat w ) | / \| X b (\hat w) \|_2^2 > 0 . $$
   If $z_j^* (\bar v) =-1 $ and $b_j (\hat w) < 0 $,  then $\hat \beta_j^0<0$ and 
  $$\beta_j^{\prime} =  - |\hat  \beta^0|  + \lambda | b_j (\hat w )| / \| X b (\hat w) \|_2^2 < 0 . $$
   If $z_j^* (\bar v) =-1 $ and $b_j (\hat w) \ge 0 $,  then $\hat \beta_j^0<0$ and 
  $$\beta_j^{\prime} =  - | \hat \beta^0|  - \lambda | b_j (\hat w ) |/ \| X b (\hat w) \|_2^2 < 0 . $$
   \fbox{Case 2} Let now $j \notin S_0$. Then
  $$\beta_j^{\prime} = - \lambda b_j (\hat w) / \| X b (\hat w) \|_2^2 , $$
  so
  $$ z_j (\hat w) \beta_j^{\prime} =  - \lambda z_j (\hat w) b_j (\hat w) / \| X b (\hat w) \|_2^2 >0 . $$
  Thus
  $$ z( \hat w) \in \partial \| \beta^{\prime} \|_1 . $$
  Furthermore, by the first part of Lemma \ref{linear.lemma}, 
  $$ X^T X  (\beta^{\prime}- \hat \beta^0 )= - \lambda X^T X  b (\hat w) / \| X b (\hat w) \|_2 =satisfies
 - \lambda  \hat W z (\hat w) .$$
 So $\beta^{\prime} =: \hat \beta$ satisfies the KKT conditions with $\tilde \zeta = 
 z (\hat w) $. 
 We further have
 \begin{eqnarray*}
  \| X (\hat \beta - \hat \beta^0)  \|_2^2  &= &
 \lambda^2 b^T (\hat w) \hat W z (\hat w)  / \| X b (\hat w ) \|_2^2 \\
 &=& \lambda^2 / \| X b (\hat w) \|^2 \\
 & \ge & \lambda^2 s_0 / (n \hat \kappa^2 ( 1+ \bar v , S_0) )
 \end{eqnarray*} 
 where in the last step we used the second part of Lemma \ref{linear.lemma}.
 Finally, by the triangle inequality
 \begin{eqnarray*}
  \| X ( \hat \beta -  \beta^0 )\|_2 &\ge& \| X ( \hat \beta - \hat \beta^0) \|_2 -{\bf U} (S_0)  \\
  &\ge & {\lambda \sqrt {s_0} \over\sqrt n  \hat \kappa (1+ \bar v, S_0)}  - {\bf U} (S_0) \\
  &\ge & {\lambda \sqrt {s_0} \over \sqrt n \hat \kappa (1+ \bar v, S_0)}  - \sqrt {s_0} - \sqrt {2x} .  \\
  \end{eqnarray*} 
 \hfill $\sqcup \mkern -12mu \sqcap$
 
  \subsubsection{Random part} \label{probabilistic2.section}
    In Lemma \ref{deterministic2.lemma}, 
     we showed that the conclusion (\ref{noisylowerbound.equation}) of 
     Theorem \ref{noisy.theorem}
   holds on the set ${\cal T}$. This subsection obtains 
that $\PP ( {\cal T}) \ge 1-  \exp [-t] + \exp[-x]$.
       
    \begin{lemma} \label{probabilistic2.lemma}
    It holds that
    $$ \PP ({\cal T} ) \ge 1- \exp[-t ] - \exp [-x] . $$
   \end{lemma}
   
   {\bf Proof of Lemma \ref{probabilistic2.lemma}.} 
     Apply Lemma \ref{Gauss.lemma} with $Z_j= \hat u_j / u_j$ for $j \in S_0$ and $Z_j = \hat v_j / v_j $ for $j \notin S_0$ to find that with probability at least $1- \exp[-t]$
   $$    | \hat u_j | \le \lambda \bar u_j \ \forall j \in S_0 , \ 
    |\hat v_j | \le \lambda \bar v_j  \ \forall j \notin S_0 . $$
    Furthermore, the random variable  ${\bf U}^2 (S_0) $ has a  chi-squared distribution with $s_0$ degrees of freedom.
    Lemma \ref{chi2.lemma} gives that with probability at least $1- \exp [-x]$, 
    $$ {\bf U} (S_0)  \le \sqrt {  s_0 } + \sqrt {2x}  . $$
 \hfill $\sqcup \mkern -12mu \sqcap$ 
 
 \subsubsection{Collecting the pieces}
   Combining  Lemma \ref{deterministic2.lemma} with Lemma \ref{probabilistic2.lemma} completes
   the proof of Theorem \ref{noisy.theorem}.

\subsection{Proof of Theorem \ref{hidden-oracle2.theorem}.}
The proof is along the lines of Theorem \ref{randomdesign.theorem}.

\subsubsection{Comparing the KKT conditions}
We compare the KKT conditions for the noisy Lasso with those
 for the noiseless Lasso.
 
       \begin{lemma}\label{hidden-oracle1a.lemma} It holds that
    $$ \| X({\hat \beta} - {\beta^*} )\|_2^2 +  \lambda \| \hat \beta \|_1 -
  \lambda^* \hat \beta^T z^* \le  (\hat \beta -\beta^* )^T X^T   \epsilon + (\lambda- \lambda^* ) \| \beta^* \|_1 .  $$
    \end{lemma}
    
    {\bf Proof of Lemma \ref{hidden-oracle1a.lemma}.}
 The KKT conditions (\ref{hatKKT.equation}) for $\hat \beta$
 can be written as
 $$ X^T X  ( \hat \beta - \beta^0 ) + \lambda \hat \zeta = X^T \epsilon  . $$
 where $\hat \zeta \in \partial \| \hat \beta \|_1 $.
 By the KKT conditions (\ref{KKT-noiseless.equation}) for $\beta^*$
 $$X^T X   ( \beta^* - \beta^0 ) + \lambda^* \zeta^* =0 . $$
 Hence, taking the difference
 $$X^T X ( \hat \beta - \beta^* ) + \lambda \hat \zeta -  \lambda^* \zeta^*  = X^T \epsilon  . $$
 Multiply by $(\hat \beta - \beta^*)^T$ to find
 $$\| X({\hat \beta } - {\beta^*} )\|_2^2 + \lambda ( \hat \beta - \beta^* )^T \hat \zeta - \lambda^* (\hat \beta - \beta^*)^T \zeta^*  =
  (\hat \beta -\beta^* )^T X^T   \epsilon . $$
 But
 \begin{eqnarray*}
 & \ & \lambda ( \hat \beta - \beta^* )^T \hat \zeta - \lambda^* (\hat \beta - \beta^*)^T \zeta^* \\
 &=& \lambda \| \hat \beta \|_1 - \lambda^* \hat \beta^T \zeta^* + \lambda^* \|\beta^{*}\|_1 - \lambda \beta^{*T} \hat \zeta \\
 &=&  \lambda \| \hat \beta \|_1 - \lambda^* \hat \beta^T \zeta^* + \lambda \|\beta^{*}\|_1 - \lambda \beta^{*T} \hat \zeta
 - (\lambda - \lambda^* ) \| \beta^* \|_1  \\
 & \ge & \lambda \| \hat \beta \|_1 - \lambda^*  \hat \beta^T \zeta^* - (\lambda- \lambda^* ) \|\beta^* \|_1 
 \end{eqnarray*}
 where we used that 
 $$ \| \beta^* \|_1 - \beta^{*T} \hat \zeta \ge 0 . $$
Therefore 
 $$ \| X({\hat \beta} - {\beta^*} )\|_2^2 +  \lambda \| \hat \beta \|_1 -
  \lambda^* \hat \beta^T z^* \le  (\hat \beta -\beta^* )^T X^T   \epsilon + (\lambda- \lambda^* ) \| \beta^* \|_1 .  $$
  \hfill $\sqcup \mkern -12mu \sqcap$
  
  \subsubsection{Projections} \label{projections3a.section}
  Recall the notation of Subsection \ref{projections3.section}.
    We let moreover $\hat v_{-S_0} $ be the vector
    $$\hat v_{-S}^S := ( X_{-S} {\rm A} X_{S} )^T \epsilon . $$
      As before, we denote the projection of $\epsilon $ on the space spanned by the columns of $X_{S}$ by
    $ \epsilon {\rm P} X_{S} $    and write
    $${\bf U} (S):= \| \epsilon {\rm P } X_{S} \|_2  .$$
    
      \subsubsection{Choice of $\lambda$} \label{lambda3a.section}
 Recall that we require that for some $t>0$    
  $$
     \lambda>  \| v_{-S}^S \|_{\infty} 
   \sqrt {2  (\log (2p) +t) } .
  $$

  \subsubsection{The set ${\cal T}^S$}
  Recall
  $$ \bar v^S := v_j^S \sqrt {2 ( \log (2p) +t ) } / \lambda , \ j \notin S . $$
  Let
  $${\cal T}^{S} := \{ |\hat v_j | \le \lambda \bar v_j \ \forall \ j \notin S \} 
\cap  \{ {\bf U} (S) \le \sqrt {s } + \sqrt {2x} \} .$$

\subsubsection{Deterministic part}
  
 \begin{lemma}\label{onTS.lemma} On the set ${\cal T}^S$ it holds that
  $$ \| X (\hat \beta - \beta^* ) \|_2 \le \sqrt { s} + \sqrt {2x} + (\lambda- \lambda^* )  \sqrt{ s/n} / \hat \kappa (\bar w^S , S) . $$
 \end{lemma}
 
 {\bf Proof of Lemma \ref{onTS.lemma}.}
  Since $S_* \subset S$
 $$X(\hat \beta- \beta^* )= X_S \hat b_S   + X_{-S} {\rm A} X_S \hat \beta_{-S} $$ 
 where
 $$ X_{S} \hat b_S = X_S (\hat \beta_S - \beta_S^*) + (X_{-S} {\rm P} X_S) \hat \beta_{-S} .$$
  In view of Lemma \ref{hidden-oracle1a.lemma},
 $$  \| X( \hat \beta - \beta^*)  \|_2^2 
 +  \lambda \| \hat \beta \|_1 -
  \lambda^* \beta^T z^* $$ $$ \le  \hat b_S^T X_S^T  \epsilon+
   \biggl [ X_{-S} {\rm A} X_S \hat \beta_S \biggr ]^T  \epsilon + (\lambda - \lambda^* )\| \beta^* \|_1  $$
 By the Cauchy-Schwarz inequality and since we are on ${\cal T}^S$
 $$\hat b_S^T X_S^T  \epsilon
  \le {\bf U} (S)  \| X \hat b_S  \|_2 \le ( \sqrt {s } + \sqrt {2x} ) \| X \hat b_S \|_2 \le (\sqrt { s} +
  \sqrt {2x} ) \| X ( \hat \beta - \beta^*) \|_2  $$
  where in the last inequality we used Pythagoras rule.
 Moreover, by the definition of $\hat v_{-S}^S$ and since we are on the set ${\cal T}^S$
 $$  \biggl [ X_{-S} {\rm A} X_S \hat \beta_{-S} \biggr ]^T  \epsilon= \hat \beta_{-S}^T \hat v_{-S}^S \le
 \lambda \sum_{j \notin S}  \bar v_{-S}^S |\hat \beta_j|  .$$
 On the other hand, 
 $$\lambda \| \hat \beta_{-S} \|_1 - \lambda^* \zeta_{-S}^{*T} \hat \beta_{-S} \ge \lambda \sum_{j \notin S} (1- \lambda^* |\zeta_j^*|/ \lambda ) |\hat \beta_j | 
$$
and
$$ ( \lambda - \lambda^* ) \| \beta^* \|_1 - \lambda \| \hat \beta_S \|_1 + \lambda^* z^{*T} \hat \beta_S 
\le ( \lambda - \lambda^* ) \| \hat \beta_S - \beta_S^* \|_1 . $$
 
     If $\| X ( \hat \beta - \beta^*)  \|_2 \le \sqrt s +\sqrt {2x} $ we are done.
     Suppose therefore that $\| X ( \hat \beta - \beta^* \|_2> \sqrt s+ \sqrt {2x} $. 
     Then we see that
     \begin{eqnarray*}
  & \ &  \| X ( \hat \beta - \beta^*) \|_2^2  - (\sqrt {s } + \sqrt {2x}) \| X (\hat \beta - \beta^*) \|_2 \\
  &=& \| X ( \hat \beta - \beta^*) \|_2 \biggl (\| X ( \hat \beta - \beta^*) \|_2- \sqrt  s   - \sqrt {2x} \biggr ) \\
  & > & 0 .
  \end{eqnarray*}
    But then
    $$ \lambda \sum_{j \notin S} (1- \bar v_j^S - \lambda^*
 |\zeta_j^*| / \lambda ) | \hat \beta_j |  < (\lambda - \lambda^* ) \| \hat \beta_S - \beta^* \|_1 . $$
    or
    $$ \| \hat \beta_S- \beta_S^* \|_1 - \| \bar W^S \hat \beta_{-S} \|_1  > 0 . $$
    Then
    $$ \| \hat \beta_S- \beta_S^* \|_1 - \| \bar W^S \hat \beta_{-S} \|_1 \le 
   ( \sqrt {s/n}  ) \| X (\hat \beta - \beta^* ) \|_2/ \hat \kappa ( \bar w^S , S) . $$
    We thus arrive at
\begin{eqnarray*}
 &  & \| X (\hat \beta - \beta^* )  \|_2^2  \\
&\le&  \biggl ( \sqrt {s} + \sqrt {2x} + (\lambda- \lambda^* )  \sqrt{ s/n} /\hat \kappa
(\bar w^S , S)  \biggr ) \| X (\hat \beta - \beta^*) \|_2
     \end{eqnarray*} 
     or
     $$ \| X (\hat \beta - \beta^* ) \|_2 \le \sqrt { s} + \sqrt {2x} + (\lambda- \lambda^* )  \sqrt {s/n} / \hat \kappa (\bar w^S , S) . $$
   \hfill $\sqcup \mkern -12mu \sqcap$
   
   \subsubsection{Random part}
   
   \begin{lemma}\label{probabilistic3.lemma} We have
   $$ \PP ( {\cal T}^S ) \ge 1- \exp[-t ] - \exp [-x] . $$
    \end{lemma}
    
    {\bf Proof of Lemma \ref{probabilistic3.lemma}.} 
    This follows from Lemma \ref{Gauss.lemma} and Lemma \ref{chi2.lemma}.
    \hfill $\sqcup \mkern -12mu \sqcap$
    
    \subsubsection{Finalizing the proof of Theorem \ref{hidden-oracle2.theorem}}
    Combine Lemma \ref{onTS.lemma} with Lemma \ref{probabilistic3.lemma}. 
    
    \subsection{Proof of the lemma in Section \ref{noisyTV.section}}
    
    {\bf Proof of Lemma \ref{weighted.lemma}.} 
     Write $g_i:= w_i f_i$, $i=1 , \ldots, n $ and $u_j := d_j /2$, $j=2 , \ldots , s $. 
          Then we have 
     $$ \sum_{j=1}^s | g_{d_j+1} - g_{d_j} | - \sum_{i=2}^{d_1} | g_i - g_{i-1} | 
    - \sum_{j=2}^{s-1}  \sum_{i=d_j+1}^{d_{j+1} } | g_{i} - g_{i-1} | - \sum_{i=d_s+1}^n | g_i- g_{i-1} | $$
      \begin{eqnarray*}
   & \le& { 1 \over d_1} \sum_{i=1}^{d_1} | g_i | + { 1 \over u_2} \sum_{i=d_1+1}^{d_1+u_2 }  | g_i  | \\
   &+& { 1 \over d_2 - u_2} \sum_{i=d_1+u_2+1}^{d_1+ d_2}  | g_i| +{1 \over u_3}
   \sum_{i=d_1+ d_2+1 }^{d_1 + d_2 + u_3} |g_i| \\
   & \ & \cdots \\
   & +& {1 \over d_{s-1} - u_{s-1}} \sum_{i=d_1 + \cdots + d_{s-2} + u_{s-1}+1}^{d_1 + \cdots + d_{s-1}} | g_i  | 
   + {1 \over u_{s}} \sum_{i=d_1 + \cdots + d_{s-1}+1}^{d_1 + \cdots + d_{s-1} + u_s  }   | g_i |\\
   &+& { 1 \over d_s -u_s} \sum_{i=d_1 + \cdots + d_{s-1} + u_{s}+1}^{d_1 + \cdots + d_{s}} | g_i  | +
   { 1 \over d_{s+1}  } \sum_{i=d_1 + \cdots + d_{s}+1}^{n}   | g_i |
   \end{eqnarray*}
      \begin{eqnarray*}
   & \le& \biggl ( { 1 \over d_1^2} \sum_{i=1}^{d_1} w_i^2  + { 1 \over u_2^2} \sum_{i=d_1+1}^{d_1+u_2 } w_i^2 \\
   &+& { 1 \over (d_2 - u_2)^2} \sum_{i=d_1+u_2+1}^{d_1+ d_2}  w_i^2 +{1 \over u_3^2}
   \sum_{i=d_1+ d_2+1 }^{d_1 + d_2 + u_3} w_i^2 \\
   & \ & \cdots \\
   & +& {1 \over (d_{s-1} - u_{s-1})^2} \sum_{i=d_1 + \cdots + d_{s-2} + u_{s-1}+1}^{d_1 + \cdots + d_{s-1}} w_i^2 
   + {1 \over u_{s}^2} \sum_{i=d_1 + \cdots + d_{s-1}+1}^{d_1 + \cdots + d_{s-1} + u_s  }   w_i^2\\
   &+& { 1 \over (d_s -u_s)^2 } \sum_{i=d_1 + \cdots + d_{s-1} + u_{s}+1}^{d_1 + \cdots + d_{s}} w_i^2 +
   { 1 \over d_{s+1}^2  } \sum_{i=d_1 + \cdots + d_{s}+1}^{n}   w_i^2 \biggr )^{1/2} \\
   &\ & \times  \biggl ( \sum_{i=1}^n f_i^2 \biggr )^{1/2} 
   \end{eqnarray*}
    \begin{eqnarray*}
   &\le&  \sqrt { {n \over d_1} + {n \over u_2} + {n \over d_2 - u_2 } + \cdots + {n \over d_{s-1}- u_{s-1}}+
   {n \over u_{s}}  + {n \over d_s - u_s } + {n \over d_{s+1} }   } \\
   & \ & \times \sqrt { \sum_{i=1}^n | f_i |^2 /n} \\
  & \ & \times  \| w \|_{\infty}. 
    \end{eqnarray*}
    Moreover
    \begin{eqnarray*}
    & \ & \sum_{j=1}^s w_{d_j+1} | f_{d_j+1} - f_{d_j} | - \sum_{i=2}^{d_1} w_i | f_i - f_{i-1} | \\
    &- & \sum_{j=2}^{s-1}  \sum_{i=d_j+1}^{d_{j+1} } w_i| f_{i} - f_{i-1} | - \sum_{i=d_s+1}^n w_i | f_i- f_{i-1} | \\
   & \le &
    \sum_{j=1}^s | g_{d_j+1} - g_{d_j} | - \sum_{i=2}^{d_1} | g_i - g_{i-1} | 
    - \sum_{j=2}^{s-1}  \sum_{i=d_j+1}^{d_{j+1} } | g_{i} - g_{i-1} | - \sum_{i=d_s+1}^n | g_i- g_{i-1} | \\
    & + &
    \sum_{i=2}^{n} | w_{i} - w_{i-1} | | f_{i-1} |,
    \end{eqnarray*} 
    and
    \begin{eqnarray*}
    \sum_{i=2}^{n} | w_{i} - w_{i-1} | | f_{i-1} | &\le&  \sqrt {\sum_{i=2}^n( w_i- w_{i-1})^2 }
    \sqrt {\sum_{i=2}^n f_{i-1}^2} \\
    & \le & \sqrt {\sum_{i=2}^n( w_i- w_{i-1})^2 }
    \sqrt {\sum_{i=1}^n f_{i}^2}
    \end{eqnarray*}
    Thus we conclude
     \begin{eqnarray*}
     & \ &  \sum_{j=1}^s w_{d_j+1} | f_{d_j+1} - f_{d_j} | \\
     &- & \sum_{i=2}^{d_1} w_i | f_i - f_{i-1} | 
    - \sum_{j=2}^{s-1}  \sum_{i=d_j+1}^{d_{j+1} } w_i | f_{i} - f_{i-1} | - \sum_{i=d_s+1}^n w_i | f_i- f_{i-1} | \\
    & \le & \left ( \| w \|_{\infty} \sqrt {{ 1 \over d_1}  + \sum_{j=2}^s { 4 \over d_j } + { 1 \over d_{s+1} } } 
    + \sqrt { n \sum_{i=2}^n (w_i - w_{i-1} )^2 } \right )  \sqrt {\sum_{i=1}^n f_i^2/n } .
    \end{eqnarray*}

   \hfill $\sqcup \mkern -12mu \sqcap$

   \subsection{Proof of Theorem \ref{tight.theorem}}\label{tight.section}
   This follows from Corollary \ref{concentration.corollary} combined with Theorem
   \ref{noiseless.theorem}, where in the latter we replace $\hat \Sigma := X^T X / n$
   by the population version $\Sigma_0$. This works because we 
   replaced Condition \ref{betaminnoiseless.condition} by its population
   counterpart Condition \ref{betamin.condition}.

 \section{Tools from probability theory}\label{tools.section}
 
 We first present three standard lemmas for Gaussian random variables,
 Lemmas \ref{Gauss.lemma}, \ref{chi2.lemma} and \ref{2dGauss.lemma}.
 These three lemmas are followed by a concentration of measure result
 and a result for Gaussian quadratic forms. 
   
   \begin{lemma} \label{Gauss.lemma}
   Let $Z_1 , \ldots , Z_p$ be standard normal random variables.
Then it holds for all $t>0$ that 
 $$ \PP \biggl ( \max_{1 \le j \le p} | Z_j | \ge \sqrt { 2 (\log (2p) + t) } \biggr )  \le \exp[-t] . $$
    \end{lemma}
    {\bf Proof of Lemma \ref{Gauss.lemma}.}
    For each $t>0$ 
    $$ \PP ( | Z_1 | \ge \sqrt {2t} ) \le 2 \exp[-t] . $$
    So by the union bound, for any $t>0$,
 \begin{eqnarray*}
   \PP \biggl ( \max_{1 \le j \le p} | Z_j | > \sqrt { 2 (\log (2p) + t) } \biggr )  
 &\le&
 p \PP(  | Z_1 | \ge \sqrt { 2 (\log (2p) + t) } )  \\
 & \le&  2 p \exp[- (\log (2p +t)] = \exp[-t] . 
 \end{eqnarray*}
 \hfill $\sqcup \mkern -12mu \sqcap$

     \begin{lemma} \label{chi2.lemma} 
     Let $Z:= (Z_1 , \ldots , Z_T)^T$ be a vector with i.i.d.
     standard Gaussian entries. Then
     it holds for all $x>0$ that
   $$\PP \biggl ( \| Z \|_2     \ge \sqrt {T}  +  \sqrt {2x}  \biggr ) \le \exp[-x] $$ and
   $$\PP\biggl  ( | \| Z \|_2     - \sqrt {T} | \ge  \sqrt {2x} \biggr  ) \le 2\exp[-x] .$$ 
    \end{lemma}
    
         {\bf Proof of Lemma \ref{chi2.lemma}.}
         This follows from concentration of measure 
          (\cite{borell1975brunn}, \cite{GineNickl}, Theorem 2.5.7) because
          the map $Z \mapsto \| Z \|_2$ is Lipschitz.
          Alternatively, one may apply
  Lemma 1 in \cite{laurent2000adaptive}. 
   \hfill $\sqcup \mkern -12mu \sqcap$
   
    \begin{lemma} \label{2dGauss.lemma}
    Let $(U,V) \in \R^{n \times 2}$ have
  i.i.d Gaussian rows with mean zero and covariance matrix
  $$\begin{pmatrix} \sigma_u^2 & \sigma_{uv} \cr  \sigma_{uv} & \sigma_v^2 \cr \end{pmatrix} .$$
  Then for all $t>0$, with probability at least $1-4\exp[-t]$
   $$ |U^T V - n \sigma_{uv} |  \le 3 \sigma_u \sigma_v \biggl ( \sqrt {2n t} +
     t \biggr ) . $$
     \end{lemma}
     
     {\bf Proof of Lemma \ref{2dGauss.lemma}.}
By standard arguments (see \cite{vdG2017} for tracking down
some constants) one can derive that with probability at least $1- 4\exp[-t]$
  $$ |U^T V - n \sigma_{uv} |  \le (  \sigma_u \sigma_v + 2 |\sigma_{u,v} |) \sqrt {2n t} +
   ( \sigma_u \sigma_v + 2 |\sigma_{u,v} |) t . $$
   We simplify this to: with probability at least $1- 4\exp[-t]$
   $$ |U^T V - n \sigma_{uv} |  \le 3 \sigma_u \sigma_v \biggl ( \sqrt {2n t} +
     t \biggr ) . $$
     \hfill $\sqcup \mkern -12mu \sqcap$
     
     This is the concentration of measure lemma that we use in Section \ref{randomdesign.section}.
   \begin{lemma}\label{concentration.lemma}
  For any $b \in \R^p$ and all $x >0$, we have 
  $$ \PP \biggl (   \|X (\hat \beta - b) \|_2 \ge m_b + \sqrt {2x} 
 \biggr ) \le  \exp[-x ]  $$
 and
  $$ \PP \biggl (  \biggl | \|X (\hat \beta - b) \|_2 - m_b \biggr | \ge \sqrt {2x} 
 \biggr ) \le 2 \exp[-x ]  $$
  where $m_b := \EE ( \|X (\hat \beta - b) \|_2  \vert X) $. 
  \end{lemma}
  
  {\bf Proof of Lemma \ref{concentration.lemma}.} 
 This follows from concentration of measure
  see e.g.\ \cite{borell1975brunn}, or \cite{GineNickl}, Theorem 2.5.7,  as the map $\epsilon \mapsto \|
  X (\hat \beta- b ) \|^2 $ is Lipschitz, see also \cite{wainwrightgeer}.
  \hfill $\sqcup \mkern -12mu \sqcap$

     Finally, we give a result for Gaussian quadratic forms.
     
     \begin{lemma}\label{quadraticform.lemma}
     Let $X$ have i.i.d. ${\cal N} (0, \Sigma_0)$-distributed rows and let
     $M$ be a (sequence of) constant(s) such that
     $$M^2 = o\biggl ( n /( \| \Sigma_0 \|_{\infty} \log (2p)) \biggr  ). $$
      Then, for a suitable sequence $\eta_M =o(1)$, 
     with probability tending to one
$$ \inf_{\| b \|_1 \le M\| \Sigma_0^{1/2}  b \|_2} { \| X b \|_2^2/ n   \over 
\| \Sigma_0^{1/2} b \|_2^2 } \ge (1- \eta_M)^2 . $$
\end{lemma}

{\bf Proof of Lemma \ref{quadraticform.lemma}.}
See for example Chapter 16 in \cite{vdG2016}
 and its references, or
\cite{GeerMuro2014}.
\hfill $\sqcup \mkern -12mu \sqcap$

\bibliographystyle{plainnat}
\bibliography{reference}
\end{document}